\documentclass[11pt,reqno]{amsart}

\usepackage[text={160mm,240mm},centering]{geometry}            
\usepackage{amssymb,amsmath,setspace,geometry,indentfirst,changepage,inputenc,amsthm}
\usepackage{mathrsfs,amsfonts,savesym,graphicx,bm,color}

\usepackage{graphicx}
\usepackage{amssymb}
\usepackage{dsfont}

\usepackage{lscape}

\usepackage{tikz}
\usetikzlibrary{positioning}

\usepackage{float}

\usepackage[all,2cell]{xy}    
\UseAllTwocells               

\usepackage{amsmath,amsthm}

\usepackage[mathscr]{eucal}
\usepackage[all]{xy}
\usepackage{mathrsfs}
\usepackage{hyperref}
\usepackage{color}

\newtheorem{theorem}{Theorem}[section]
\newtheorem{lemma}[theorem]{Lemma}

\newtheorem{conjecture}[theorem]{Conjecture}
\newtheorem{corollary}[theorem]{Corollary}
\newtheorem{proposition}[theorem]{Proposition}

\newtheorem{definition-lemma}[theorem]{Definition-Lemma}
\newtheorem{definition-theorem}[theorem]{Definition-Theorem}

\theoremstyle{definition}
\newtheorem{example}[theorem]{Example}
\newtheorem{definition}[theorem]{Definition}

\newtheorem{remark}[theorem]{Remark}

\usepackage{fancyhdr}
\allowdisplaybreaks

\title{On Motivic Zeta Functions and Stringy E-function via Embedded $\mathbb{Q}$-Resolution}

\author{Yifan Chen}
\address{
	Department of Mathematical Sciences,
	Tsinghua University,
	Beijing, 100084, P. R. China.}
\email{c-yf20@mails.tsinghua.edu.cn}

\author{Quan Shi}
\address{
	Department of Mathematical Sciences,
	Tsinghua University,
	Beijing, 100084, P. R. China.}
\email{shiq24@mails.tsinghua.edu.cn / thusq20@gmail.com}

\author{Huaiqing Zuo}
\address{Department of Mathematical Sciences,
	Tsinghua University,
	Beijing, 100084, P. R. China.}
\email{hqzuo@mail.tsinghua.edu.cn}

\begin{document}
	
	\maketitle
	
	%
	
	\begin{abstract}
		We provide the formula of motivic zeta function for semi-quasihomogeneous singularities and in dimension two, we determine the poles of zeta functions. We also give another formula for stringy E-function using embedded $\mathbb{Q}$-resolution, and we utilize it to calculate the stringy E-function for semi-quasihomogeneous polynomials and non-degenerate polynomials. 
	\end{abstract}
	
	\tableofcontents
	
	\subsection*{Notation and Convention} (1) We use $\bm x$ to mean the array of indeterminants $(x_1,...,x_n)$. (2) For $a \in \mathbb N_{>0}$, $\zeta_a := e^{2\pi\sqrt{-1}/a}$ is the unit root. (3) we use $\frac{1}{a}(a_1,...,a_n)$ to mean the cyclic group generated by $\mathrm{diag}\{\zeta_a^{a_1},...,\zeta_{a}^{a_n}\}$ and the same notation to mean the cyclic quotient singularity regarding this group. $\mathrm{diag}\{\zeta_a^{a_1},...,\zeta_{a}^{a_n}\}$ is called the generator of this group. (4) For a function $f : X \to \mathbb C$, with $X$ a variety, set $Z(f) := \{f = 0\}$. (5) $\chi(X)$ is Euler characteristic of a variety $X$. (6) For a weight $\bm w\in \mathbb N_{>0}^n$, we always require $\mathrm{gcd}(w_1,...,w_n) = 1$.

	\section{Introduction}\label{sec1}
	
	For $f\in \mathbb C[\bm x]\setminus \mathbb C$, or more generally an effective divisor on a smooth complex variety, there is a well-known conjecture, the monodromy conjecture, which predicts the poles of motivic zeta function are always in the logarithmic set of eigenvalues of local monodromies, or even roots of Bernstein-Sato polynomials (\cite{Motivic_Igusa_Zeta_Function}). The monodromy conjecture is one of the main open problems in singularity theory, relating to many areas of mathematics, including number theory, motivic integration, birational geometry, $\mathscr D$-module theory, etc.
	
	The first ingredient of the monodromy conjecture is the motivic zeta funciton. It is an analogue of Igusa's $p$-adic zeta function (\cite{Igusa_Zeta_Function_Prologue}), added by Kontsevich's idea of motivic integration. It is defined as below.
	\begin{definition}
		For a subvariety $W$ of $\mathbb C^n$, the motivic zeta function of $f$ at $W$ is 
		\begin{displaymath}
			Z_{f,W}^{\mathrm{mot}}(s) := \int_{\rho^{-1}(W)} (\mathbb L^{-s})^{\mathrm{ord}_t f(\bm u)} \,\mathrm{d}\mu(\bm u),
		\end{displaymath}
		where $\mathrm{ord}_t f(\bm u)$ means that for $\bm u \in J_{\infty}(X)$, the order of the power series $f(\bm u)$, and $\rho:J_{\infty}(\mathbb{C}^{n}) \rightarrow \mathbb{C}^{n}$ is the projection. For the meaning of the integration notation, we refer to Section \ref{sec_Grothendieck_Ring} for details.
	\end{definition} 
	
	In \cite{Motivic_Igusa_Zeta_Function}, Denef and Loeser proved $Z_{f,W}^{\mathrm{mot}}(s)$ is ratinoal in $\mathbb L^{-s}$.
	\begin{theorem}[Denef-Loeser, \cite{Motivic_Igusa_Zeta_Function}]\label{Ch2_Sec3_Motivic_and_Log_Resolution_Thm}
		Let $\pi : Y\to \mathbb C^{n}$ be a log resolution of $Z(f)\subseteq \mathbb C^n$ with numerical data $(N_i,\nu_i)_{i=0}^r$ i.e. $\mathrm{div}(f \circ \pi) = \sum_{i=0}^r N_iE_i$ and $K_\pi = \sum_{i=0}^r (\nu_i-1)E_i$ being SNC divisors, then 
		\begin{displaymath}
			Z_{f,W}^{\mathrm{mot}}(s) = \mathbb L^{-n}\sum_{I\subseteq \{0,1,...,r\}} [E_I^\circ\cap \pi^{-1}(W)] \prod_{i\in I} \frac{(\mathbb L-1)\mathbb L^{-N_is-\nu_i}}{1-\mathbb L^{-N_is-\nu_i}}, \tag{\ref{Ch2_Sec3_Motivic_and_Log_Resolution_Thm}.1}
		\end{displaymath}
		where $E_I^\circ := (\bigcap_{i\in I} E_i) \setminus (\bigcup_{j\notin I} E_j)$. In particular, $Z_{f,W}^{\mathrm{mot}}(s)$ is a rational function.
	\end{theorem}
	\begin{remark}
		$Z_{f,W}^{\mathrm{mot}}(s)$ can be specialized to the topological zeta function by
		\begin{displaymath}
			Z_{f,W}^{\mathrm{top}}(s) = \sum_{I\subseteq \{0,1,...,r\}} \chi(E_I^\circ\cap \pi^{-1}(W)) \prod_{i\in I} \frac{1}{N_is+\nu_i}.
		\end{displaymath}
	\end{remark}
	One can define the poles of $Z_{f,W}^{\mathrm{mot}}$ accordingly. Although we obtain a set of candidate poles from (\ref{Ch2_Sec3_Motivic_and_Log_Resolution_Thm}.1), many of them are cancelled in pratice. Such cancellation is very mysterious and interesting (\cite{Motivic_Zeta_Function_on_Q_Gorenstein_Varieties_Leon_Martin_Veys_Viu_Sos}). It is one of the main difficulties in the study of monodromy conjecture and so it is of great value to find poles that survive.
	
	Another two ingredients, local monodromies (\cite{Milnor_Fibration_Monodromy_2}) and Bernstein-Sato polynomials (\cite{Rationality_of_Roots_of_B-Function_Kashiwara, SZ24a}) are important invariants from the Milnor fibration and $\mathscr D$-module theory. They are both very difficult to compute, too. 
	\begin{conjecture}[Monodromy Conjecture, \cite{Motivic_Igusa_Zeta_Function}]
		Let $f\in \mathbb C[\bm x]\setminus \mathbb C$. If $s_0$ is a pole of $Z_{f,W}^{\mathrm{mot}}(s)$, then (1) $(s-s_0)\mid b_f(s)$; (2) $e^{s_0\cdot 2\pi \sqrt{-1}}$ is an eigenvalue of a local monodromy at some neighborhoods of points in $Z(f)\cap W$.
	\end{conjecture}
	Historically, the monodromy conjecture evolved from the fact that poles of the archimedean zeta function is always a root of the Bernstein-Sato polynomial (\cite{Introduction_to_the_Monodromy_Conjecture_Veys_2024, SZ24b}). Then Igusa developed an analogue, the $p$-adic zeta function (\cite{Igusa_Zeta_Function_Prologue}) and found a similar phenomenon. Then Denef and Loeser specialized Igusa's zeta function to a $p$-adic one and posed a topological version (\cite{Topoloical_Zeta_Function_Independent_of_Resolution}). Both were unified by a motivic version in \cite{Motivic_Igusa_Zeta_Function}. However, nothing is known in the literature how motivic zeta function is related to these two invariants in general. The monodromy conjecture then turns out rather miraculous.
	
	\vspace{0.5em}
	
	Therefore, it is of great value to know more about the motivic zeta function. For example, compute an explicit formula. In this paper, we will treat the case when $f$ is semi-quasihomogeneous. The definition is as follows:
	\begin{definition}[Quasihomogeneous and semi-quasihomogeneous polynomial]
		A polynomial $f \in \mathbb{C}[x_{1},...,x_{n}]$ is called quasihomogeneous if there exists a weight $\bm{w}=(w_{1},...,w_{n}) \in \mathbb{N}^{n}$ such that for each monomial $x_{1}^{\alpha_{1}}...x_{n}^{\alpha_{n}}$ in $f$, $\alpha_{1}w_{1}+...+\alpha_{n}w_{n}=d$ for some fixed $d>0$. In this case, we say $f$ is quasihomogeneous of degree $d$ with weight $\bm{w}$. $f$ is called semi-quasihomogeneous if $f=f_{d}+g$, where $f_{d}$ is quasihomogeneous with weight $\bm{w}$ of degree $d$ and $f_{d}$ defines an isolated singularity at $\bm{0}$, and $g$ is of degree bigger than $d$ under the weight $\bm{w}$. Such $\bm{w}$ is called the weight of $f$ and $f_{d}$ is called the principal part of $f$.
	\end{definition}
	
	Semi-quasihomogeneous singularities form a very large type of singularities and abundant information can be read from the weight. Restricting to the case $f$ is quasihomogeneous, the topological type of $(Z(f),\bm 0) \subseteq (\mathbb C^n,\bm 0)$ is closely related by its weight type (see \cite{XY89}). Fix a family $\mathfrak Q \subseteq \mathbb C\{\bm x\}$ of isolated singularities at $\bm 0$. An analytic invariant $\sigma(g)$ for $g\in \mathbb C\{\bm x\}$ is called a topological invariant of $\mathfrak Q$ if $\sigma(g_1) = \sigma(g_2)$ as long as $(V(g_1),\bm 0)$ is homeomorphic to $(V(g_2),\bm 0)$ for $g_{1}, g_{2} \in \mathfrak{Q}$. 
	
	\vspace{0.5em}

	Suppose $f \in \mathbb C[\bm x]\setminus \mathbb C$ is a semi-quasihomogeneous polynomial with a unique singularity at the origin. Suppose $\bm w$ is the weight of $f$ and $\deg_{\bm w} f = d$. Let $\Pi : \widehat{\mathbb{C}^{n}}(w) \longrightarrow \mathbb{C}^{n}$ be the weighted blow-up of $\mathbb{C}^{n}$ (also see \cite{An_Introduction_to_p_adic_and_Motivic_Integration_Vius_Sos}). More precisely,
	
	\begin{displaymath}
		\widehat{\mathbb{C}^{n}}(\bm w)=\{(\bm x,[\bm u]_{\bm w}) \in \mathbb{C}^{n} \times \mathbb{P}^{n-1}_{\bm w}| \ \exists \lambda \in \mathbb C^*,x_{i}=\lambda^{w_{i}}u_{i}, \forall i\},
	\end{displaymath}
	and $\Pi$ is the the projection to $\bm x$. $\Pi$ is an isomorphism outide $\bm 0$ by definition. Moreover, $\widehat{\mathbb C^n}(\bm w)$ is naturally a $V$-manifold, i.e. locally an abelian quotient singularity (see \cite{Mixed_Hodge_on_the_Vanishing_Cohomology}). For $\widehat{\mathbb C^n}(\bm w)$ the covering $U_i = \{(\bm x,\bm u) \mid u_i \neq 0\},i=1,...,n$ satisfies $U_i \simeq \mathbb C^n/G_i$, $G_i = \frac{1}{w_i}(w_1,...,w_{i-1},-1,w_{i+1},...,w_n)$. This quotient map $\pi_i : \mathbb C^n \to \mathbb C^n/G_i$ is identified with
	\begin{align*}
		&(u_1,...,u_{i-1},\lambda_{i},u_{i+1},...,u_n) \\
		&\mapsto  (\lambda_{i}^{w_1}u_1,...,\lambda_{i}^{w_{i-1}}u_{i-1},\lambda_{i}^{w_i},\lambda_{i}^{w_{i+1}}u_{i+1},...,\lambda_{i}^{w_n}u_n,[u_1,...,u_{i-1},1,u_{i+1},...,u_n]).
	\end{align*}
	Taking $x_{j}=\lambda_i^{w_{j}}u_{j},j \neq i$ and $x_i = \lambda_i^{w_i}$ into the expression of $f$, we get 
	\begin{displaymath}
		f(\bm x) =\lambda_i^{d}(f_{d}(u_{1},...,u_{i-1},1,u_{i+1},...,u_{n})+\lambda_i h(\lambda,u_{1},...,u_{i-1},1,u_{i+1},...,u_{n})) ,
	\end{displaymath}
	where $f_d$ is the principal part of $f$. Let $\hat H$ be the strict transform of $\{f = 0\}$ via $\Pi$. Then the inverse image of $\hat{H}\cap U_i$ under $\pi_i$ is exactly $\{H_i := f(\bm x)\cdot \lambda_i^{-d} = 0\}$. 
	

	For $\mathbb C^n/G_i$, we find a natural stratification $\bigsqcup_{I\subseteq \{1,....n\}} \mathbb C_I/G_{i}$ with $G_I*\mathbb C_I = \mathbb C_I \subseteq \mathbb C^n$ and  the stablizer of each point in $\mathbb C_I $ is the same, denoted as $G_{I,i} \subseteq G_i$, where $\mathbb{C}_I = A_I^1 \times ... \times A_I^n \subseteq \mathbb C^n$, with $A_I^i = \mathbb{C}^{*}$ if $i \in I$, and $A_I^i = \{0\}$ if $i \notin I$(see \textbf{Subsection \ref{Subsec_Quotient_Singularity}}). Let $U_{I,i}$ be the image of $\mathbb C_I/G_i$ in $U_i$, and $E$ be the exceptional divisor. Suppose the inverse image of $E \cap U_i$ and $\hat H \cap U_i$ under $\mathbb C^n \to \mathbb C^n/G_i$ are $E_i'$ and $\hat H_i'$. We further set $P_{I,i} : = \mathbb C_I\setminus (H_i'\cup E_i')$, $\hat{H}_{I,i}^\circ := (\mathbb C_I \cap H_i') \setminus E_i'$, $E_{I,i}^\circ = (E_i' \setminus \hat H_i')\cap \mathbb C_I$,  and $\hat{E}_{I,i} = \mathbb C_I \cap H_i'\cap E_i'$. Let $Z_i = \{i+1,...,n\}$.  For the chart $U_{i}$ and $I \subset Z_{i}$, we denote the set of points $p$ in $U_{i} \cap \{\lambda=f_{d}=0\}$ such that $\frac{\partial f_{d}}{\partial u_{j}} \neq 0, \frac{\partial f_{d}}{\partial u_{l}}=0$ ($l,j \neq i, 1 \le l \le j-1$) $E_{i,j}$.
	
	We now state the Theorem A of this paper. Let $G_i' = \frac{1}{w_i}(-1,w_1,...,\hat w_i,...,w_n) \simeq G_i$, with generators corresponded. Take $G_{I,i}'$ the image of $G_{I,i}$ under this isomorphism.  Let $G_{i,j}'' = \frac{1}{w_i}(-1,d,w_1,...,\hat w_i,...,\hat w_j,...,w_n)$ be another group isomorphic to $G_i$, with generators corresponded as well. Let $G_{I,i,j}$ be the image of $G_{I\setminus \{j\},i}$ under this isomorphism.
	
	
	
	
	We also use the following notation as in \cite{Motivic_Zeta_Function_on_Q_Gorenstein_Varieties_Leon_Martin_Veys_Viu_Sos}. For a finite group $G$ consisting of diagonal matrices and $\gamma \in G_{k}$, $\gamma$ can be written uniquely as $\mathrm{diag}(\zeta_{r}^{\varepsilon_{\gamma,1}},...,\zeta_{r}^{\varepsilon_{\gamma,n}})$, where $r=|G|$, $\zeta_{r}=e^{\frac{2\pi \sqrt{-1}}{r}}$ and $0 \le \varepsilon_{\gamma,i} \le r-1$. For any $\bm{k}=(k_{1},...,k_{n}) \in \mathbb{Q}^{n}$, let $\varpi_{\bm{k}}: G \longrightarrow \mathbb{Q}$ given by $\varpi_{\bm{k}}(\gamma)=\frac{1}{r}\sum_{i=1}^{n}k_{i}\varepsilon_{\gamma,i}$ and  $S_{G}(\bm{N},\bm{\nu},s):=\sum_{\gamma \in G}\mathbb{L}^{\varpi_{\bm N}(\gamma)s+\varpi_{\bm \nu}(\gamma)}$.
	
	So far we can formula the motivic zeta function.
	
	\vspace{0.5em}
	\noindent\textbf{Theorem A }(Theorem \ref{motivic zeta function})\textbf{.} 
	\textit{Suppose $f \in \mathbb C[\bm x]\setminus \mathbb C$ is a semi-quasihomogeneous polynomial with a unique singularity at the origin. With notations as above, the motivic zeta function of $f$ is given as}
	\begin{align*}
		\mathbb L^n\cdot Z_{f}^{\mathrm{mot}}(s) & = \mathbb L^{n}-[Z(f)] + ([Z(f)]-1)\cdot \frac{(\mathbb L-1)\mathbb L^{-(s+1)}}{1-\mathbb L^{-(s+1)}} \\
		& \quad + \sum_{i=1}^n\sum_{I \subseteq Z_i} ((\mathbb L-1)^{\vert I\vert}-[\hat E_{I,i} / G_{i}]) \cdot S_{G_{I,i}'}(\bm N_3, \bm \nu_3,s) \cdot \frac{(\mathbb L-1)\mathbb L^{-(ds+\vert \bm w\vert)}}{1-\mathbb L^{-(ds+\vert \bm w\vert)}}\\
		& \quad + \sum_{i=1}^n\sum_{I \subseteq Z_i} [ (E_{i,j} \cap \mathbb{C}_{I}) / G_{i}] \cdot S_{G_{I,i,j}}(\bm N_4, \bm \nu_4,s) \cdot \frac{(\mathbb L-1)\mathbb L^{-(s+1)}}{1-\mathbb L^{-(s+1)}}\cdot \frac{(\mathbb L-1)\mathbb L^{-(ds+\vert \bm w\vert)}}{1-\mathbb L^{-(ds+\vert \bm w\vert)}}.
	\end{align*}
	\textit{Here $\bm N_3,\bm N_4,\bm \nu_3, \bm \nu_4$ are from \textbf{Table \ref{Tab_Local_Q_Normal_Crossings}}. The local motivic zeta function at $\bm 0$ is given as}
	\begin{align*}
		\mathbb L^n\cdot Z_{f,\bm 0}^{\mathrm{mot}}(s) & = \sum_{i=1}^n\sum_{I \subseteq Z_i} ((\mathbb L-1)^{\vert I\vert}-[\hat E_{I,i} / G_{i}]) \cdot S_{G_{I,i}'}(\bm N_3, \bm \nu_3,s) \cdot \frac{(\mathbb L-1)\mathbb L^{-(ds+\vert \bm w\vert)}}{1-\mathbb L^{-(ds+\vert \bm w\vert)}} \\
		& \quad + \sum_{i=1}^n\sum_{I \subseteq Z_i} [(E_{i,j} \cap \mathbb{C}_{I}) / G_{i}] \cdot S_{G_{I,i,j}}(\bm N_4, \bm \nu_4,s) \cdot \frac{(\mathbb L-1)\mathbb L^{-(s+1)}}{1-\mathbb L^{-(s+1)}}\cdot \frac{(\mathbb L-1)\mathbb L^{-(ds+\vert \bm w\vert)}}{1-\mathbb L^{-(ds+\vert \bm w\vert)}}.
	\end{align*}
	\textit{Moreover, $Z_{f,\bm 0}^{\mathrm{mot}}(s) = Z_{f_d,\bm 0}^{\mathrm{mot}}(s)$, where $f_d$ is the principal part of $f$.}
	\vspace{0.5em}
	
	Applying the formula in dimension two, we obtain an interesting result.
	
	\vspace{0.5em}
	\noindent\textbf{Theorem B }(Corollary \ref{dimension 2 result} and \ref{dimension 2 non-cancellation})\textbf{.} 
	\textit{Suppose $f \in \mathbb C[x_1,x_2]$ is semi-quasihomogeneous and has a unique singularity at $\bm 0$. Let $\bm w = (p,q)$, $(p,q) = 1$, be the weight of $f$ and $d = \deg_{\bm w} f$. Suppose $f(x_1^{p},x_1^{q}u) = x_1^d f_1(u)$ and $f(x_2^{p}v,x_2^{q}) =x_2^{d}f_2(v)$, then}
	\begin{align*}
		&\mathbb L^2\cdot Z_{f,\bm 0}^{\mathrm{mot}}(s) \\
		& = \frac{(\mathbb L-1)\mathbb L^{-(ds+\vert \bm w\vert)}}{1-\mathbb L^{-(ds+\vert \bm w\vert)}}\cdot 
		\begin{cases}
			\mathbb L-1-\frac{d}{pq}+S_1+S_2 + \frac{d}{pq}\frac{(\mathbb L-1)\mathbb L^{-(s+1)}}{1-\mathbb L^{-(s+1)}}, & f_1(0)\neq 0, f_2(0) \neq 0.\\
			\mathbb L-1-\frac{d-p}{pq}+S_2 + (\frac{d-p}{pq}+S_1')\frac{(\mathbb L-1)\mathbb L^{-(s+1)}}{1-\mathbb L^{-(s+1)}}, & f_1(0) = 0, f_2(0) \neq 0.\\
			\mathbb L-1-\frac{d-q}{pq}+S_1 + (\frac{d-q}{pq}+S_2')\frac{(\mathbb L-1)\mathbb L^{-(s+1)}}{1-\mathbb L^{-(s+1)}}, & f_1(0) \neq 0, f_2(0) = 0.\\
			\mathbb L-1-\frac{d-p-q}{pq}+(\frac{d-p-q}{pq}+S_1'+S_2')\frac{(\mathbb L-1)\mathbb L^{-(s+1)}}{1-\mathbb L^{-(s+1)}}, & f_1(0) = 0, f_2(0) = 0.
		\end{cases}	
	\end{align*}
	\textit{Here $S_i = S_{G_i}(\bm N_3,\bm \nu_3,s),\ S_i' = S_{G_i'}(\bm N_4,\bm \nu_4,s)$, and $G_1 = \frac{1}{p}(-1,q),\ G_1 = \frac{1}{q}(-1,p),\ G_1' = \frac{1}{p}(-1,d),\ G_2' = \frac{1}{q}(-1,d).$ In particular, $Z_{f,\bm 0}^{\mathrm{mot}}(s)$  depends only on the weight type and $\bm w$-degree of $f$.}
	
	\noindent\textit{Moreover, suppose $d > p+q$, then $-1$ and $=\vert \bm w\vert/d$ are both poles of $Z_{f,\bm 0}^{\mathrm{mot}}(s)$.}
	\vspace{0.5em}
	
	The above theorem tells us that $Z_{f,\bm 0}^{\mathrm{mot}}(s)$ depends only on $(\bm w,d)$ if $f$ is quasihomogenesous in $2$ variables. It is very surprising since in general even $Z_{f,\bm 0}^{\mathrm{top}}(s)$ is not a toplogical invariant of the family of polynomials (see \cite{ABCN$^+$02}). In higher dimension cases, it seems not so likely that $Z_{f,\bm 0}^{\mathrm{mot}}(s)$ is a topological invariant for quasihomogeneous singularities, but $Z_{f,\bm 0}^{\mathrm{top}}(s)$ may be. Combining \textbf{Example \ref{homogeneous motivic zeta function}}, we propose a conjecture as below.
	\begin{conjecture}\label{topological zeta function topological}
		Let $\mathbf{QH}_n$ be the family of all quasihomogeneous polynomials in $\mathbb C[\bm x]$ with a unique singularity at $\bm 0$. For $f\in \mathbf{QH}_n$ with weight type $\bm w$ and $d = \deg_{\bm w} f$. Then 
		
		(1) $Z_{f,\bm 0}^{\mathrm{top}}(s)$ depends only on $(\bm w,d)$.
		
		(2) $Z_{f,\bm 0}^{\mathrm{top}}(s)$ is a topological invariant of the family $\mathbf{QH}_n$.
	\end{conjecture}
	
	By means of motivic integration, Batyrev introduced the ``stringy E-function" as an invariant of singularity in \cite{Batyrev_Stringy_Invariant} and \cite{Stringy_Hodge_Number}. This invariant is defined only on log terminal singularities and Batyrev utilized them to test topological mirror symmetry for singular Calabi-Yau varieties. He also gave the definition of stringy hodge numbers and proved a version of Macky correspondence using stringy E-function. The definition is as follows.
	
	\begin{definition}[Stringy E-function]
		Suppose $X$ is a log terminal variety of dimension $n$. Take a log resolution $\phi: Y \longrightarrow X$, and let $D_{i}(i \in S)$ be its irreducible components of exceptional divisor and $a_{i}(i \in S)$ be its log discrepancies. For $I \subset S$, we denote $\mathring{D_{I}}:=\bigcap_{i \in I}D_{i} \backslash \bigcup_{i \notin I}D_{i}$, then the stringy E-function, i.e. the stringy E-function of $X$ is defined to be:
		\begin{flalign*}
			E_{st}(X):=\sum_{I \subset S}E(D_{I})\prod_{i \in I}\frac{uv-1}{(uv)^{a_{i}}-1},
		\end{flalign*}
		where $E(\cdot)$ means the E-function of Hodge.
		Using Theorem $2.15$ in \cite{Introduction_to_Motivic_Integration}, we can rewrite the stringy E-function of $X$ into the following integration:
		\begin{flalign*}
			E_{st}(X)=E(\mathbb{L}^{n}\int_{J_{\infty}(Y)}\mathbb{L}^{-\mathrm{ord}_{t}K_{\phi}}d\mu_{J_{\infty}(Y)}),
		\end{flalign*}
		where $K_{\phi}=K_{Y}-\phi^{*}(K_{X})$ is the relative canonical divisor.
	\end{definition}
	In this paper, we will use embedded $\mathbb{Q}$-resolution to calculate stringy E-function, with the following theorem:
	
	\vspace{0.5em}
	\noindent\textbf{Theorem C }(Theorem \ref{stringy E-function for Q-Gor})\textbf{.} 
	\textit{	Assume $X$ is a $\mathbb{Q}$-Gorenstein variety with at worst log terminal singularity of pure dimension $n$. Let $\pi : Y \longrightarrow X$ be an embedded $\mathbb{Q}$-resolution such that $Y=\bigsqcup_{k \ge 0} Y_{k}$, and for any point $q \in Y_{k}$, we have $(Y,q) \simeq (\mathbb{C}^{n}/G_{k}, 0)$, and the relative canonical divisor $K_{\pi}$ is locally given by $x_{1}^{\nu_{1,k}-1}...x_{n}^{\nu_{n,k}-1}$, where $x_{1},...,x_{n}$ are the coordinates of $\mathbb{C}^{n}$. Then the stringy E-function of $X$ is given by:}
	\begin{flalign*}
		E_{st}(X)=\sum_{k \ge 0}E(Y_{k})E(S_{G_{k}}(\bm{\nu}_{k}))\prod_{i=1}^{n}\frac{uv-1}{(uv)^{\nu_{i,k}}-1},
	\end{flalign*}
	\textit{where the notation $S_{G_k}(\nu_k) := S_{G_k}(\bm 0,\bm \nu_k,0)$.}
	
	As the applications of the formula, we can use it to calculate the stringy E-function of semi-quasihomogeneous polynomials and non-degenerate polynomials. Here are two theorems about them:

	\vspace{0.5em}
	\noindent\textbf{Theorem D }(Theorem \ref{stringy E-function for semi-quasihomogeneous})\textbf{.} 
	\textit{	Let $f \in \mathbb{C}[x_{1},...,x_{n}]$ a semi-quasihomogeneous polynomial with a unique singularity at $0$ and $H:=\{f=0\}$. The notations $\hat{H},U_{i},f_{d},u_{i}$ are the same as in Theorem A. We denote the part of $(U_{i}-\bigcup_{l=1}^{i-1}U_{l}) \cap \hat{H}$ where $\frac{\partial f_{d}}{\partial u_{j}} \neq 0, \frac{\partial f_{d}}{\partial u_{1}}=...=\frac{\partial f_{d}}{\partial u_{j-1}}=0$ by $E_{i,j}$. For the point $q \in E_{i,j}$, we use $\lambda,u_{1},...,\hat{u_{i}},...,\hat{u_{j}},...,u_{n}$ as local coordinates of $(\hat{H},q)$, where $\hat{u_{j}},\hat{u_{i}}$ means $u_{j},u_{i}$ are omitted. For $I \subset \{1,...,\hat{i},...,\hat{j},...,n\}$, when $q \in E_{i,j} \cap \mathbb{C}_{I}$, we have $(\hat{H},q) \cong (\mathbb{C}^{n-1}/G_{i,I,j},0)$ for some $G_{i,I,j}$. Then the stringy E-function of $H$ is given by:
		\begin{flalign*}
			E(H)-1+\sum_{i=1}^{n}\sum_{j \neq i}\sum_{I \subset \{1,...,\hat{i},...,\hat{j},...,n\} }E((E_{i,j} \cap \mathbb{C}_{I})/G_{i})E(S_{G_{i,I,j}}(|\bm w|,1,...,1))\frac{uv-1}{(uv)^{|\bm w|}-1}.
		\end{flalign*}
	}

	Before we state the theorem about the stringy E-function of non-degenerate polynomial, we will introduce some notations. Let $f \in \mathbb{C}[x_{1},...,x_{n}]$ be a non-degenerate polynomial with a unique singularity at $0$ and $H=\{f=0\} \subset \mathbb{C}^{n}$. Let $\Sigma(f)'$ be a simplicial subdivision of the normal fan of $f$ without changing the set of one-dimensional cones and $\Pi: X_{\Sigma(f)'} \longrightarrow \mathbb{C}^{n}$ is the corresponding morphism of toric varieties. We denote the strict transform of $H$ by $\hat{H}$ and $\pi:=\Pi|_{\hat{H}}:\hat{H} \longrightarrow H$. Suppose $\Sigma(f)'=\bigcup_{i=1}^{s}\sigma_{s}$, where $\sigma_{1},...,\sigma_{s}$ are cones of dimension $n$, and we assume the $1$-dim cones of $\sigma_{i}$ is given by $\rho_{i,1},...,\rho_{i,n}$. Then $X_{\Sigma(f)'}=\bigcup_{i=1}^{s}U_{\sigma_{i}}=\bigsqcup_{i=1}^{s}(U_{\sigma_{i}}-U_{\sigma_{1}}\cup...\cup U_{\sigma_{i-1}})$, where $U_{\sigma_{i}}$ is the affine toric variety associated to $\sigma_{i}$ and $U_{\sigma_{i}} \simeq \mathbb{C}^{n}/G_{\sigma_{i}}$ for some group $G_{\sigma_{i}}$. Let $V_{i}:=U_{\sigma_{i}}-U_{\sigma_{1}}\cup...\cup U_{\sigma_{i-1}}$ be a closed subset of $U_{\sigma_{i}}$, and $E_{i,j}$ be the subset of $V_{i} \cap \hat{H}$ such that $\frac{\partial g_{\sigma_{i}}}{\partial u_{j}} \neq 0$ and $\frac{\partial g_{\sigma_{i}}}{\partial u_{1}}=...=\frac{\partial g_{\sigma_{i}}}{\partial u_{j-1}}=0$, where $g_{\sigma_{i}}$ is the defining equation for $\hat{H}$ on $U_{i}$. For $I \subset \{1,...,\hat{j},...,n\}$ and point $p \in \mathbb{C}_{I}\cap E_{i,j}$, we have $(\hat{H},p)\simeq (\mathbb{C}^{n-1}/G_{i,I,j},0)$ for some $G_{i,I,j}$.
	
	\vspace{0.5em}
	\noindent\textbf{Theorem E }(Theorem \ref{stringy E-function for nondegenerate})\textbf{.} 
	\textit{	With the notations as above, the stringy E-function of $H$ is given by:
		\begin{flalign*}
			\sum_{i=1}^{s}\sum_{j=1}^{n}\sum_{I \subset \{1,...\hat{j},...,n\}}E((E_{i,j}\cap \mathbb{C}_{I})/G_{\sigma_{i}})S_{G_{i,I,j}}(\bm{\nu}_{i,I,j})\prod_{l \in I}\frac{uv-1}{(uv)^{\rho_{i,l}\cdot (1,...,1)-\phi_{f}(\rho_{i,l})}-1},
		\end{flalign*}
		where the $l$-th coordinate of $\bm{\nu}_{i,I,j}$ is $\rho_{l}\cdot (1,...,1)-\phi_{f}(\rho_{l})$ if $l \in I$, and $1$ otherwise.
	}

	\subsection*{Acknowledgement}
	
	We thank Nero Budur for discussions. Zuo is supported by NSFC Grant 12271280.

	\section{Preliminary}\label{sec2}
	
	\subsection{Jet Scheme, Grothendieck Ring and Motivic Integration}\label{sec_Grothendieck_Ring}
	In this section, we will introduce some basic facts about jet scheme, Grothendieck ring and motivic integration.
	\begin{definition}[Jet Scheme]
		Let $X$ be a scheme over a field $k$, for every $m \in \mathbb{N}$, consider the functor from $k$-schemes to set
		\begin{flalign*}
			Z \mapsto \mathrm{Hom}(Z \times_{\mathrm{Spec}(k)} \mathrm{Spec} (k[t]/(t^{m+1})),X).
		\end{flalign*}
		There is a $k$-scheme represents this functor, which is called the $m$-th jet scheme of $X$, denoted by $J_{m}(X)$(see Theorem 2.1, \cite{Jet_Schemes_Ishii}), i.e.
		\begin{flalign*}
			\mathrm{Hom}(Z,J_{m}(X))=\mathrm{Hom}(Z \times_{\mathrm{Spec}(k)} \mathrm{Spec} (k[t]/(t^{m+1})),X).
		\end{flalign*}
	\end{definition}
	
	For $1 \le i \le j$, the truncation map $k[t]/(t^{j}) \to k[t]/(t^{i})$ induces natural projections between jet schemes $\psi_{i,j}: J_{j}(X) \to J_{i}(X)$. If we identify $J_{0}(X)$ with $X$, we can get natural projection $\pi_{i}: J_{i}(X) \to X$. One can check that $\{J_{m}(X)\}_{m}$ is an inverse system and the inverse limit $J_{\infty}(X):= \varprojlim_m J_{m}(X)$ is called the arc space of $X$. There are also natural projections $\psi_{i}: J_{\infty}(X) \to J_{i}(X)$. If $X$ is of finite type, $J_{m}(X)$ is also of finite type for each $m \in \mathbb{N}$, but $J_{\infty}(X)$ is usually not.
	
	The Grothendieck ring is a feasible way to make the set of all algebraic varieties an abelian group.
	\begin{definition}[Grothendieck Ring]\label{Ch2_Sec2_Grothendieck_Ring}
		Let $k$ be a field and let $\mathrm{Var}_{k}$ be the category of $k$-varieties. The Grothendieck group of $k$-varieties, $K_0[\mathrm{Var}_k]$, is defined to be the quotient group of the free abelian group with basis $\{[X]\}_{X\in\mathrm{Var}_k}$, modulo the following relations.
		\begin{align*}
			& [X]-[Y], \ X \simeq Y\\
			& [X]-[X_{\mathrm{red}}]\\
			& [X]-[U]-[X\setminus U],\ U\subseteq X\ \mathrm{open}.
		\end{align*}
		One can further define a multiplication structure on $K_0[\mathrm{Var}_k]$ by
		\begin{displaymath}
			[X] \cdot [Y] := [X\times Y].
		\end{displaymath}
		This makes $K_0[\mathrm{Var}_k]$ a ring, called the Grothendieck ring of $k$-varieties.
		
		Let $\mathbb L = [\mathbb A_k^1]$ and $K_0[\mathrm{Var}_k]_{\mathbb L}$ be the localization at $\mathbb L$.
	\end{definition}
	
	\noindent\textbf{Kontsevich's Completion}
	
	We use $\mathcal M$ to mean $K_0[\mathrm{Var}_{k}]$ and $\widetilde{\mathcal M}$ to mean $K_0[\mathrm{Var}_{k}]_{\mathbb L}$. There is a natural decreasing filtration $F^\bullet$ on $\widetilde{ \mathcal M}$. For $m\in \mathbb Z$, $F^m$ is the subgroup of $\widetilde{M}$ generated by $[S] \cdot \mathbb L^{-i}$ with $\dim S- i \leq -m$. It is actually a ring filtration i.e. $F^m \cdot F^n \subseteq F^{m+n}$.
	\begin{definition}
		The Kontsevich's completed Grothendieck ring $\widehat{\mathcal M}$ is defined to be 
		\begin{displaymath}
			\widehat{\mathcal M} := \varprojlim_{m\in \mathbb Z} \widetilde{\mathcal M}/F^m.
		\end{displaymath}
	\end{definition}
	
	Now we fix an algebraic variety $X$ of pure dimension $d$ and we will do the motivic integration on the space $J_{\infty}(X)$. Firstly we need to define the measurable sets, which turns out to be the cylinders in $J_{\infty}(X)$, with the following definition:
	\begin{definition}[Cylinder]
		A subset $C \subset J_{\infty}(X)$ is called a cylinder if $C=\psi_{m}^{-1}(A)$ for some constructible subset $A \subset J_{m}(X)$. 
	\end{definition}
	
	We want to define the measure $\mu$ on $J_{\infty}(X)$ such that $\mu(C)=\frac{[\psi_{m}(C)]}{\mathbb{L}^{md}}$ for a cylinder $C=\psi_{m}^{-1}(A)$. To show that this definition is independent of the choice of $m$, we use the following theorem:
	
	\begin{theorem}[\cite{Motivic_Integration_on_Arbitrary_Varieties_Denef_Loeser}]
		If $C \subset J_{\infty}(X)$ is a cylinder and $C \cap J_{\infty}(X_{\mathrm{sing}})=\emptyset$, then there exists integer $m$ such that:     
		
		(1) $\psi_{m}(C)$ is constructible and $C=\psi_{m}^{-1}(\psi_{m}(C))$.
		
		(2) For $n \ge m$, the projection $\psi_{n+1}(C) \to \psi_{n}(C)$ is a piecewise trivial fibration with fiber $\mathbb{A}^{d}$.
		
		In particular, $\frac{[\psi_{n}(C)]}{\mathbb{L}^{nd}}$ are equal for $n \ge m$.
	\end{theorem} 
	
	For general cylinder $C$, the condition of the above theorem may not hold and we will define the measure with value in Kontsevich's completed Grothendieck ring $\widehat{\mathcal M}$, with the following theorem:
	
	\begin{theorem}[\cite{Motivic_Integration_on_Arbitrary_Varieties_Denef_Loeser}]
		Let $C$ be a cylinder of $J_{\infty}(X)$, then the limit 
		\begin{flalign*}
			\mu(C):=\lim_{m \rightarrow \infty} \frac{[\psi_{m}(C)]}{\mathbb{L}^{md}}
		\end{flalign*}
		exists in $\widehat{\mathcal M}$.
	\end{theorem}
	
	Now we can define the motivic integration:
	\begin{definition}
		Suppose $C \subset J_{\infty}(X)$ is a cylinder and $\alpha: C \to \mathbb{Z} \cup \{\infty\}$ with cylinder fiber, we define the motivic integration
		\begin{flalign*}
			\int_{C}\mathbb{L}^{-\alpha}d\mu=\sum_{n \in \mathbb{Z}}\mu(\alpha^{-1}(n))\mathbb{L}^{-n},
		\end{flalign*}
		if the right hand side convergences.
	\end{definition}
	
	\begin{remark}
		In \cite{Motivic_Zeta_Function_on_Q_Gorenstein_Varieties_Leon_Martin_Veys_Viu_Sos}, the authors gave another measure $\mu^{\mathbb{Q}-Gor}$ on $J_{\infty}(X)$ for those $X$ which are $\mathbb{Q}$-Gorenstein, and they also define a motivic zeta function using their new measure. It can be shown that when $X$ is smooth, these two measures coincide and we will see in the stringy E-function part that $\mu^{\mathbb{Q}-Gor}$ is the "right" measure on $\mathbb{Q}$-Gorenstein varieties. 
	\end{remark}

	\subsection{(Embedded) $\mathbb Q$-Resolution and Computation Method}
	
	As in \textbf{Theorem \ref{Ch2_Sec3_Motivic_and_Log_Resolution_Thm}}, motivic zeta function can be computed via log resolutions. However, in some cases we do not have to resolve the singularity completely. Things also work if there are some quotient singularities. Hence we introduce the notion of embedded $\mathbb{Q}$-resolution, which is a main tool for this paper.
	
	\begin{definition}
		Let $X$ be a $\mathbb{Q}$-Gorenstein variety of pure dimension $n$ with at worst log terminal singularities, $H \subset X$ is a subvariety of codimension one, an embedded $\mathbb Q$-resolution of $H$ is a proper morphism $\mu : Y \to X$ with the following.
		
		(1) $Y$ has only abelian quotient singularities.
		
		(2) $\mu$ is an isomorphism over $Y\setminus \mu^{-1}(H_{\mathrm{sing}})$.
		
		(3) There is a natural stratification of $Y = \bigsqcup_{k = 1}^m Y_k$ such that for each $p\in Y$, $(Y, p) \simeq (\mathbb C^n/G_k,0)$ for some abelian group, and the coordinates $x_{1},...,x_{n}$ on $\mathbb C^n/G_k$ satisfying locally the relative canonical divisor $K_{\mu}$ is given by $x_{1}^{\nu_{k,1}-1}...x_{n}^{\nu_{k,n}-1}$, $\mu^{*}(H)$ is given by $x_{1}^{N_{k,1}}...x_{n}^{N_{k,n}}$ and $G_{k}$ acts on $\mathbb{C}^{n}$ diagonally. Here $ N_{k,i}, \nu_{k,i} \in \mathbb Q_{\geq 0}$ and $G_k$  depend only on $k$. 
	\end{definition}
	
	\begin{remark}
		The definition above also works when we take $H$ to be $0$, in that case, we say $Y$ is an embedded $\mathbb{Q}$-resolution of $X$.
	\end{remark}
	
	\begin{remark}
		In the next subsection, we see that weighted blow up is always an embedded $\mathbb Q$-resolution of semi-quasihomogeneous singularities.
	\end{remark}
	
	When we consider a polynomial $f \in \mathbb{C}[\bm x]$, the motivic zeta function of $f$ can be computed via embedded $\mathbb Q$-resolution, by taking $X$ to be $\mathbb{C}^{n}$ and $H=\{f=0\}$.  The main theorem of \cite{Motivic_Zeta_Function_on_Q_Gorenstein_Varieties_Leon_Martin_Veys_Viu_Sos} provides a formula of $\mathbb{Q}$-Gorenstein version of motivic zeta function for two divisors on $\mathbb Q$-Gorenstein varieties. We only state case when $D_1 = \mathrm{div}\, f, D_2 = 0, X = \mathbb C^n$. Note that in this case this motivic zeta function coincides with the one we defined since $\mathbb{C}^{n}$ is smooth.
	
	\begin{theorem}[\cite{Motivic_Zeta_Function_on_Q_Gorenstein_Varieties_Leon_Martin_Veys_Viu_Sos}, $D_1 = \mathrm{div}\, f, D_2 = 0, X = \mathbb C^n$]\label{Thm_Formula_of_Zeta_Function_via_Embedded_Q_Resolution}
		Suppose $\mu : Y \to \mathbb C^n$ is an embedded $\mathbb Q$-resolution of $f$. $Y_k,G_k,\bm N_k,\bm \nu_k\in \mathbb Q_{\geq 0}$ are information described in the above paragraphs. Then for every subvariety $W\subseteq \mathbb C^n$, we have
		\begin{displaymath}
			Z_{f,W}^{\mathrm{mot}}(s) = \mathbb L^{-n} \sum_{k=1}^m [Y_k \cap \mu^{-1}(W)] \cdot S_{G_k}(\bm N_k,\bm \nu_k,s) \prod_{i=1}^n \frac{(\mathbb L-1)\mathbb L^{-N_{k,i}s-\nu_{k,i}}}{1-\mathbb L^{-N_{k,i}s-\nu_{k,i}}} \in \widehat{\mathcal M}[[\mathbb L^{-s/R}]].
		\end{displaymath}
		Here $R$ is a large positive integer such that $R\bm N_k,R\bm \nu_k \in \mathbb Z^n$ for all $k$. 
	\end{theorem}

	\section{Zeta Function for Semi-Quasihomogeneous Singularities}

	\subsection{A Natural Stratification for Abelian Quotient Singularities}\label{Subsec_Quotient_Singularity}
	
	Suppose $G\subseteq \mathrm{GL}_n(\mathbb C)$ is a finite abelian group. Since matrices in $G$ commute with each other, they can be diagonalized simultaneously. Hence, we may assume $G$ is a group of diagonal matrices. 
	
	Consider $\mathbb C^n/G$. There is a natural $\mathbb G_m^n = (\mathbb C^*)^n$-action on $\mathbb C^n/G$, simply by 
	\begin{displaymath}
		(\lambda_1,...,\lambda_n) * [(z_1,...,z_n)] = [(\lambda_1z_1,...,\lambda_nz_n)].
	\end{displaymath}
	Hence, there is no difference for points on each $\mathbb G_m^n$-orbit. Since $G\subseteq \mathbb G_m^n$, $G$-orbits are always $\mathbb G_m^n$-orbits. Consider a natural stratification on $\mathbb C^n$ i.e. 
	\begin{displaymath}
		\mathbb C^n = \bigsqcup_{I \subseteq \{1,...,n\}} \mathbb C_I.
	\end{displaymath}
	Here, $\mathbb{C}_I = A_I^1 \times ... \times A_I^n \subseteq \mathbb C^n$, with $A_I^i = \begin{cases}
		\{0\}, & i\not\in I\\
		\mathbb C^*, & i\in I
	\end{cases}$. This stratification is exactly the decomposition of $\mathbb C^n$ into $\mathbb G_m^n$-orbits. Therefore, 
	\begin{displaymath}
		\mathbb C^n/G = \bigsqcup_{I \subseteq \{1,...,n\}} \mathbb C_I/G
	\end{displaymath}
	is also a stratification. The following proposition tells us geometrically that for points of $[\bm z] \in \mathbb C_I/G$, $(\mathbb C^n/G,[\bm z])$ are also quotient singularities of a group associated with $I$.
	\begin{proposition}\label{Prop_Other_Points_of_Quotient_Singularity_are_again_Quotient_Singularity}
		For an arbitrary point $[\bm z] \in \mathbb C_I/G$, there is an analytic neighborhood $U\subseteq \mathbb C_I/G$ isomorphic to a neighborhood of $[\bm 0]_{G_{[\bm z]}} \in \mathbb C^n/G_{[\bm z]}$, with $[\bm z]$ mapped to the origin. Here, $G_{[\bm z]}\subseteq G$ is the stablizer of $[\bm z]$.
	\end{proposition} 
	\begin{proof}
		Since there is a $\mathbb G_m^n$-action on $\mathbb C^n/G$, we may assume coordinates of $\bm z$ are either $0$ or $1$. Furthermore, after swapping coordinates, assume $\bm z = (1,1,..,1,0,...,0)$, with $1$ takes $r$ places. Clearly, $G_1 := G_{[\bm z]}$ consists of those matrices with $1$ filled on the first $r$ positions of the diagonal. 
		
		Let $\tilde U = \mathbb B_{\varepsilon}(\bm z) \subseteq \mathbb C^n, \varepsilon \ll 1$ be a neighborhood of $\bm z$ and $U$ its image in $\mathbb C^n/G$. Let $V$ be the image of $\mathbb B_{\varepsilon}(\bm 0) \subseteq \mathbb C^n$ in $\mathbb C^n/G_{1}$. To avoid confusion, we write the corresponding group after the bracket for representatives. Consider the morphism
		\begin{displaymath}
			\rho : U \to V, [(1+v_1,...,1+v_r,v_{r+1},...,v_{n})]_G \mapsto [(v_1,...,v_n)]_{G_1}.
		\end{displaymath}
		One can easily write an inverse of $\rho$ in the same formula, since we have picked $\varepsilon$ small enough so that no two different copies $g\mathbb B_{\varepsilon}(\bm z),g\in G$ intersect with each other.
	\end{proof}

	Next, we prove that for an abelian finite subgroup $G\subseteq \mathrm{GL}_n(\mathbb C)$, one has $[\mathbb C^n/G] = \mathbb L^n \in \mathcal M$. We think it may be familiar to experts, but it is hard to find a complete reference in the literature. So we give a  brief proof here.
	
	\begin{lemma}\label{Cyclic_Singularity_Grothendieck_Class}
		Suppose $G = \langle g\rangle$ is a finite cyclic subgroup of $\mathrm{GL}_n(\mathbb C),n\geq 0$, then $[(\mathbb C^*)^n/G] = (\mathbb L-1)^n$. Note that $[\{pt\}/G] = 1$, consequently, $[\mathbb C^n/G] = \mathbb L^n$.
	\end{lemma}
	\begin{proof}
		The case $n = 0$ is trivial. We may assume $g$ is diagonal and $n \geq 1$. Consider the natural stratification of $\mathbb C^n$ under $\mathbb G_m^n$ action i.e. $\mathbb C^n = \bigsqcup_{I\subseteq \{1,...,n\}} \mathbb C_I$, where $\mathbb C_I$ is a product of some $\{0\}$s and $\mathbb{C}^{*}$s, with $i$-th component $C^*$ for all $i\in I$. Suppose the stablizer for elements in $\mathbb{C}_I$ are $G_I$. It suffices to prove $[\mathbb C_I/G_I] = [\mathbb C_I]$ and it is exactly the first assertion since $G_I$ is also cyclic.
		
		Without loss of generality, we consider $X = \mathbb C_{\{1,...,n\}}$. Suppose $g = \mathrm{diag}\{\zeta^{a_1},...,\zeta^{a_n}\}$, where $\zeta$ is a $d$-th primitive unit root and $\mathrm{gcd}(a_1,...,a_n) = 1$. Suppose $a_1k - a_2l = (a_1,a_2)$, and take $a_1' = a_1/(a_1,a_2),a_2' = a_2/(a_1,a_2)$. Pick $b_1,b_2,c_1,c_2\in \mathbb Z$ such that
		\begin{displaymath}
			\begin{pmatrix}
				k & l \\
				a_2' & a_1'
			\end{pmatrix}
			\begin{pmatrix}
				b_1 & b_2\\
				c_1 & c_2
			\end{pmatrix}
			= 
			\begin{pmatrix}
				b_1 & b_2\\
				c_1 & c_2
			\end{pmatrix}
			\begin{pmatrix}
				k & l \\
				a_2' & a_1'
			\end{pmatrix}
			= I_2.
		\end{displaymath}
		Here we use $a_1'k-a_2'l = 1$ and hence coefficients of the inverse matrix are also integers. We have
		\begin{displaymath}
			\begin{pmatrix}
				b_1 & b_2\\
				c_1 & c_2
			\end{pmatrix}
			= 
			\begin{pmatrix}
				a_1' & -l\\
				-a_2' & k
			\end{pmatrix}.
		\end{displaymath}
		
		Take $g' = \mathrm{diag}\{\zeta^{(a_1,a_2)},1,\zeta_3^{a_3},...,\zeta^{a_n}\}$, we have an isomorphism $X/\langle g\rangle \overset{\simeq}{\longrightarrow} X/\langle g' \rangle$ and its inverse as below.
		\begin{align*}
			& [\lambda_1,...,\lambda_n] \mapsto [\lambda_1^{k}\lambda_2^{-l},\lambda_1^{a_2'}\lambda_2^{-a_1'},\lambda_3,...,\lambda_n].\\
			& [\mu_1,...,\mu_n] \mapsto [\mu_1^{{b_1}}\mu_2^{b_2},\mu_1^{-c_1}\mu_2^{-c_2},\mu_3,...,\mu_n].
		\end{align*}
		One can easily check they are well-defined and inverse to each other. Therefore, by induction we have $X/\langle g\rangle \simeq X/\langle g'\rangle$, where $g' = \mathrm{diag}\{\zeta,1,...,1\}$. It is $\mathrm{Spec}\, [x_1^d,x_2,...,x_n]_{x_1^dx_2...x_n}$, isomorphic to $\mathbb C_I$.
	\end{proof}

	\begin{proposition}\label{Quotient_Singularity_Grothendieck_Class}
		For any finite abelian group $G\subseteq \mathrm{GL}_n(\mathbb C)$, $[(\mathbb C^*)^{n+1}/G] = (\mathbb L-1)^{n+1}$ and $[\mathbb C^n/G] = \mathbb L^n$ for all $n\geq 0$.
	\end{proposition}
	\begin{proof}
		An finitely generated abelian group can be decomposed into direct sums of cyclic group. Suppose $G = \bigoplus_{j = 1}^r \langle g_j \rangle$, where $g_j$ are diagonal matrices. By \textbf{Lemma \ref{Cyclic_Singularity_Grothendieck_Class}}, $(\mathbb C^*)^n/\langle g_1\rangle$ is isomorphic to $(\mathbb C^*)^n$. Moreover, from the procedure of the lemma, we see $G/\langle g_1\rangle$ also acts on $(C^*)^n/\langle g_1\rangle$ diagonally. By induction, we are done.   
	\end{proof}
	
	\subsection{Formula of Zeta Function}\label{subse_zeta_function}
	
	Let $f \in \mathbb{C}[\bm x]$ be a semi-quasihomogeneous polynomial and $H:=\{f=0\} \subset \mathbb{C}^{n}$. By the definition of semi-quasihomogeneous polynomial, $f=f_{d}+g$, where $f_{d}$ is a quasihomogeneous polynomial of degree $d$ with weight $\bm w \in \mathbb N_{>0}^n$. $f_d$ defines an isolated singularity at $\bm 0$ and the degree of $g$ with respect to this weight is larger than $d$. We know that $\bm 0$ is also the unique singularity of $f$.
	
	Let $E$ be the exceptional locus of the weighted blow up $\Pi : \widehat{\mathbb C^n}(\bm w) \to \mathbb C^n$, $\hat{H}$ be the strict transform of $H$, and $\pi:=\Pi|_{\hat{H}}:\hat{H} \longrightarrow H$ be the restriction of $\Pi$ into $\hat{H}$, then $\hat{E}:=E\cap \hat{H}$ is the exceptional locus of $\pi$. When we take $x_{i}=\lambda^{w_{i}}u_{i}$ into the expression of $f$, we can get $f(\lambda^{w_{1}}u_{1},...,\lambda^{w_{n}}u_{n})=\lambda^{d}(f_{d}(u_{1},...,u_{n})+\lambda h(\lambda,u_{1},...,u_{n}))$, where $h(\lambda,u_{1},...,u_{n})=\lambda^{-d-1}g(\lambda^{w_{1}}u_{1},...,\lambda^{w_{n}}u_{n})$, so $\hat{H}=\{(\bm x),[\bm u]_{\bm w}) \in \widehat{\mathbb{C}^{n}}(\bm w)|f_{d}+\lambda h=0\}$.

	The following is a refined result of Theorem 1.2 of \cite{Monodromy_Conjecture_for_Semi-Quasihomogeneous_Hypersurfaces_Blanco_Budur_vdV}.
	\begin{theorem}\label{Blw_is_Embeddedd_Q_Resolution}
		$\Pi$ is an embedded $\mathbb Q$-resolution of $(\mathbb C^n,f)$ and $\pi$ is a $\mathbb Q$-resolution of $H$.
	\end{theorem}
	\begin{proof}
		Let $U_1,...,U_n$ be standard charts of $Y$, where 
		\begin{displaymath}
			U_i = \{(\bm{X}_i,[\bm{U_i}]_{\bm w})\mid \lambda_{i},u_1,..,\hat{u_i},...,u_n \in \mathbb C\},
		\end{displaymath} 
		where we denote $(\lambda_{i}^{w_1}u_1,...,\lambda_{i}^{w_{i-1}}u_{i-1},\lambda_{i}^{w_i},\lambda_{i}^{w_{i+1}}u_{i+1},...,\lambda_{i}^{w_n}u_n)$ as $\bm{X}_i$, $(u_1,...,u_{i-1},1,u_{i+1},...,u_n)$ as $\bm{U}_i$. Besides, we denote $(u_1,...,u_{i-1},\lambda_{i},u_{i+1},...,u_n)$ as $\bm{U}_i^x$.
		
		$U_i$ is isomorphic to $\frac{1}{w_i}(w_1,...,-1,...,w_n) = \mathbb C^n/\langle \mathrm{diag}\{\zeta_{w_i}^{w_1},...,\zeta_{w_i}^{-1},...,\zeta_{w_i}^{w_n}\}\rangle$, with $i$-th position $-1$. The quotient map is given by
		\begin{align*}
			\phi_i &\, : \mathbb C^n \to \frac{1}{w_i}(w_1,...,-1,...,w_n), \bm{U}_i^x \mapsto (\bm{X}_i,[\bm{U}_i]_{\bm w}).
		\end{align*}
		One can easily see $\Pi^{-1}(\bm 0) \cap U_i = \hat E\cap U_i = \{x_i = 0\}$ is the exceptional divisor.

		On the chart $U_i$, we have
		\begin{displaymath}
			(\phi_i \circ \Pi|_{U_i})^* f = f(\bm{X}_i) = \lambda_{i}^d\cdot (f_d(\bm{U}_i) + \lambda_{i} h(\bm{U}_i^x)),
		\end{displaymath}
		where $h \in \mathbb C[\bm U_i^x]$ is a polynomial. Let $g_i = \mathrm{diag}\{\zeta_{w_i}^{w_1},...,\zeta_{w_i}^{-1},...,\zeta_{w_i}^{w_n}\}$ and $G_i =  \langle g_i\rangle \subseteq \mathrm{GL}_n(\mathbb C)$. An important observation is that $g_i$ acts on $f(\bm{X}_i)$ trivially. In fact,
		\begin{displaymath}
			f(\bm{X}_i) * g_i = f((\zeta_{w_i}^{-1}\lambda_{i})^{w_1}(\zeta_{w_i}^{w_1}u_1),...,(\zeta_{w_i}^{-1}\lambda_{i})^{w_i},...,(\zeta_{w_i}^{-1}\lambda_{i})^{w_n}(\zeta_{w_i}^{w_n}u_n)) = f(\bm {X}_i).
		\end{displaymath}
		Since $g_i$ acts on $\lambda_{i}$ as $\lambda_{i} \mapsto \zeta_{w_i}^{-1}$$\lambda_{i}$, $g_i$ acts on $H_i := f_d(\bm{U}_i) + \lambda_{i} h(\bm{U}_i^x)$ as $H_i * g_i = \zeta_{w_i}^d \cdot H_i$.
		
		Next, we show $\{H_i=0\}$ is smooth on $\mathbb C^n$. 
		
		$\{f_d(\bm{U}_i)=0\}$ is smooth on $\mathbb C^{n-1}$. We write $f_d = f_d(\bm{U_i}^x)$. Since $f_d$ is weighted homogeneous, we have
		\begin{displaymath}
			f_d = (\sum_{j\neq i} w_ju_j \frac{\partial f_d}{\partial u_j} (\bm{U}_i)) + w_i \frac{\partial f_d}{\partial x_i} (\bm{U}_i).
		\end{displaymath} 
		Hence, the equation $f_d(\bm{U_i}) = \frac{\partial f_d}{\partial u_j} (\bm{U}_i) = 0, j\neq i$ is equivalent to $\frac{\partial f_d}{\partial x_i} (\bm{U}_i) = \frac{\partial f_d}{\partial u_j} (\bm{U}_i) = 0$. There is no solution since $f_d$ has a unique singularity at the origin. Hence, $H_i$ is smooth on $\{x_i = 0\} = \hat E\cap U_i$. Outside $\hat E$, since $\Pi: \widehat{\mathbb C^n}(\bm w)\setminus \hat E \overset{\simeq}{\longrightarrow} \mathbb C^n\setminus\{\bm 0\}$ is an isomorphism and $f$ is smooth on $\mathbb C^n\setminus \{\bm 0\}$, $H_i$ is also smooth.

		The relative canonical divisor can be computed as
		\begin{displaymath}
			(\phi_i \circ \Pi|_{U_i})^*(\mathrm{d} \bm{y}) = w_i \lambda_{i}^{\vert \bm w\vert-1} \mathrm{d}\bm{U}_{i}^x\ \Rightarrow K_{\Pi}|_{U_i} = (\vert \bm w\vert-1)\hat E.
		\end{displaymath}
		
		To prove $\Pi$ is an embedded $\mathbb Q$-resolution, we have to show: (1) for each point $q \in U_i$, there is a neighborhood $V\subseteq U_i$ s.t. $(V,q) \simeq (\mathbb C^n/G_q,[\bm 0]_{G_q})$ for some abelian group $G_q \subseteq \mathrm{GL}_n(\mathbb C)$; (2) there are coordinates $\bm t$ of $\mathbb C^n$ s.t. $G_q$ acts diagonally on $\bm t$ and $\Pi^* f = \bm{t}^{\bm N}, K_\Pi = \bm t^{\bm{\nu}-\bm{e}}$, where $\bm N,\bm{\nu} \in \mathbb Q_{\geq 0}$ and $\bm e = (1,..,1)$ ($\mathbb Q$-normal crossing). 
		
		Now we apply the stratification in \textbf{subsection \ref{Subsec_Quotient_Singularity}}, say for $I \subseteq \{1,...,n\}$, $\mathbb C_I$ is the set of points in $\mathbb C^n$ with all $j$-th coordinate non-zero, $j\in I$, and others zero. Let $G_{I,i}$ be the stablizer of $\bm z\in \mathbb C_I$ in $G_i$, consisting of those matrices in $G_i$ s.t. $1$ is on all $j$-th positions of diagonal, $j\in I$ ($G_{I,i}$ depends only on $I$). By \textbf{Proposition \ref{Prop_Other_Points_of_Quotient_Singularity_are_again_Quotient_Singularity}}, $(\mathbb C^n/G_i,[\bm z]_{G_i}) \simeq (\mathbb C^n/G_{I,i},[\bm 0]_{G_{I,i}})$. Since $\{f_d(\bm{U}_i) = 0\}$ is smooth on $\mathbb C^{n-1}$, for all $\bm z \in \{H_i = 0\}$, the linear term of $H_i$ at $\bm z$ is linearly independent with $\lambda_{i}-z_i$. Then one can pick $n-2$ coordinates among $\{u_1,...,\hat{u_i},...,u_n\}$ so that they, together with $\lambda_{i}-z_i$, $H_i$, satisfies (1), (2).
		
		More precisely, on $\mathbb C_I/G_i$, the $\bm N$ and $\bm{\nu}$ are given in \textbf{Table \ref{Tab_Local_Q_Normal_Crossings}}. We take $\lambda_{i}-z_i$ and $H_i$ the first two coordinate. One has to slight revise the group $G_{I,i}$ since the coordinate change happens. The corresponding groups under the new coordinates are exactly $G_{I,i}'$ and $G_{I,i,j}$ defined in the introduction part.
		
		Since $H_i$ is smooth and $E = \{\lambda_{i} = 0\}$ intersects $\{H_i = f_d + \lambda_{i} h = 0\}$ transversally, we see that $\pi$ is a $\mathbb Q$-resolution.
		%
	\end{proof}
	
	Changing to the notations in the introduction, we have the following.
	\begin{align*}
		& P_{I,i} = \mathbb C_I \setminus (\{H_i = 0\}\cup \{\lambda_{i} = 0\}), & \hat{H}_{I,i}^\circ = \mathbb C_I\cap (\{H_i = 0\}\setminus \{\lambda_{i} = 0\}), \\
		& E_{I,i}^\circ = \mathbb C_I\cap (\{\lambda_{i} = 0\}\setminus \{H_i = 0\}), &\hat E_{I,i} = \mathbb C_I\cap \{H_i = 0\} \cap \{\lambda_{i} = 0\}.\ \ 
	\end{align*}
	The information of the embedded $\mathbb Q$-resolution are given as below.
	\begin{table}[h!]
		\centering
		\begin{tabular}{|c|c|c|}
			\hline
			
			$\bm z \in \mathbb C_I$ & $\bm N$ & $\bm{\nu}$\\
			
			\hline
			
			$P_{I,i}$ & $(0,...,0)$ & $(1,...,1)$\\
			
			\hline
			
			$H_{I,i}^\circ$ & $(0,1,0,...,0)$ & $(1,1,...,1)$\\
			
			\hline
			
			$E_{I,i}^\circ$ & $(d,0,...,0)$ & $(\vert \bm w\vert,1,...,1)$\\
			
			\hline 
			
			$\hat E_{I,i}$ & $(d,1,...,0)$ & $(\vert \bm w\vert,1,...,1)$\\
			\hline
		\end{tabular}
		\vspace{0.5em}
		\caption{Local $\mathbb Q$-normal Crossings}
		\label{Tab_Local_Q_Normal_Crossings}
	\end{table}

	Before computing the zeta function, we can slight reduce the number of strata.
	
	\begin{lemma}\label{slightly_reduction}
		Suppose $g = \mathrm{diag}\{\zeta_a^{-1},\zeta_a^{b_2},...,\zeta_a^{b_n}\} \in \mathrm{GL}_n(\mathbb C), a,b_2,...,b_n \in \mathbb Z_{>0}, b_i \leq a$, and $G = \langle g\rangle$. Then $\mathbb C^n/G$ can be decomposed as $(A/G) \sqcup \bigsqcup_{I \subseteq \{2,...,n\}} (\mathbb C_I/G)$, where $A = \mathbb C^* \times \mathbb C^{n-1}$. 
		
		The stablizer of points in of $A$ is trivial. For $\mathbb C_I$, let $d_I = \mathrm{gcd}_{j\in I} b_j$ and $b_I = a/(a,d_I)$, then the stablizer of points in $\mathbb C_I$ is $\langle g^{b_I}\rangle$. If $I = \varnothing$, it is $G$.
	\end{lemma} 
	\begin{proof}
		Let $a \in \mathbb C_I, I\subseteq \{2,...,n\}$, then $g^b$ stablizes $a$ if and only if $a \mid b_ib$ for all $i\in I$. It is equivalent to $a \mid d_I \cdot b$, or $b_I \mid b$. 
	\end{proof}
	
	
	Next we give the formula for $Z_{f}^{\mathrm{mot}}(s)$. To make the statement shorter, some notations may be provided before.  As in the proof of \textbf{Theorem \ref{Blw_is_Embeddedd_Q_Resolution}}, let $g_i = \mathrm{diag}\{\zeta_{w_i}^{w_1},...,\zeta_{w_i}^{-1},...,\zeta_{w_i}^{w_n}\}$ and $G_i = \langle g_i\rangle$. Let $Z_i = \{i+1,...,n\}$ and $A_i = \{0\}^{i-1} \times \mathbb C^* \times \mathbb C^{n-i}$. For $I \subseteq Z_i$, let $d_{I,i} = \mathrm{gcd}_{j\in I} w_j$ and $b_{I,i} = w_i/(w_i,d_{I,i})$ and $G_{I,i} = \langle g_i^{b_{I,i}} \rangle$. For the chart $U_{i}$ and $I \subset Z_{i}$, we denote the points $p$ in $\{H_{i}=0\}$ such that $\frac{\partial H_{i}}{\partial u_{j}} \neq 0, \frac{\partial H_{i}}{\partial u_{l}}=0$($l,j \neq i, 1 \le l \le j-1$) $E_{i,j}$, and when $p \in \mathbb{C}_{I}$, assume $(\mathbb{C}^{n}/G_i,p) \simeq (\mathbb{C}^{n}/G_{I,i,j},0)$. Let $\bm N_1,...,\bm N_4$ and $\bm \nu_1,...,\bm \nu_4$ be the four types of $\bm N,\bm \nu$ in \textbf{Table \ref{Tab_Local_Q_Normal_Crossings}} on the four rows respectively.
	\begin{theorem}\label{motivic zeta function}
		The motivic zeta function of $f$ is given as
		\begin{align*}
			\mathbb L^n\cdot Z_{f}^{\mathrm{mot}}(s) & = \mathbb L^{n}-[Z(f)] + ([Z(f)]-1)\cdot \frac{(\mathbb L-1)\mathbb L^{-(s+1)}}{1-\mathbb L^{-(s+1)}} \\
			& \quad + \sum_{i=1}^n\sum_{I \subseteq Z_i} ((\mathbb L-1)^{\vert I\vert}-[\hat E_{I,i} / G_{i}]) \cdot S_{G_{I,i}'}(\bm N_3, \bm \nu_3,s) \cdot \frac{(\mathbb L-1)\mathbb L^{-(ds+\vert \bm w\vert)}}{1-\mathbb L^{-(ds+\vert \bm w\vert)}}\\
			& \quad + \sum_{i=1}^n\sum_{I \subseteq Z_i} [(\mathbb C_I\cap E_{i,j}) / G_{i}] \cdot S_{G_{I,i,j}}(\bm N_4, \bm \nu_4,s) \cdot \frac{(\mathbb L-1)\mathbb L^{-(s+1)}}{1-\mathbb L^{-(s+1)}}\cdot \frac{(\mathbb L-1)\mathbb L^{-(ds+\vert \bm w\vert)}}{1-\mathbb L^{-(ds+\vert \bm w\vert)}}.
		\end{align*}
		And the local motivic zeta function at $\bm 0$ is given as
		\begin{align*}
			\mathbb L^n\cdot Z_{f,\bm 0}^{\mathrm{mot}}(s) & = \sum_{i=1}^n\sum_{I \subseteq Z_i} ((\mathbb L-1)^{\vert I\vert}-[\hat E_{I,i} / G_{i}]) \cdot S_{G_{I,i}'}(\bm N_3, \bm \nu_3,s) \cdot \frac{(\mathbb L-1)\mathbb L^{-(ds+\vert \bm w\vert)}}{1-\mathbb L^{-(ds+\vert \bm w\vert)}} \\
			& \quad + \sum_{i=1}^n\sum_{I \subseteq Z_i} [(\mathbb C_I\cap E_{i,j}) / G_{i}] \cdot S_{G_{I,i,j}}(\bm N_4, \bm \nu_4,s) \cdot \frac{(\mathbb L-1)\mathbb L^{-(s+1)}}{1-\mathbb L^{-(s+1)}}\cdot \frac{(\mathbb L-1)\mathbb L^{-(ds+\vert \bm w\vert)}}{1-\mathbb L^{-(ds+\vert \bm w\vert)}}.
		\end{align*}
		Moreover, $Z_{f,\bm 0}^{\mathrm{mot}}(s) = Z_{f_d,\bm 0}^{\mathrm{mot}}(s)$, where $f_d$ is the principal part of $f$.
	\end{theorem}
	\begin{proof}
		We have the following stratification of $\widehat{\mathbb C^n}(\bm w)$.
		\begin{displaymath}
			\widehat{\mathbb C^n}(\bm w) = \bigsqcup_{i=1}^n \bigg((A_i/G_i)\sqcup \bigsqcup_{I \subseteq Z_i'} (\mathbb C_I/G_{I,i})\bigg).
		\end{displaymath}
		Adopt the notations in \textbf{Theorem \ref{Blw_is_Embeddedd_Q_Resolution}}. On the chart $U_i$ with fake coordinate $\lambda_i,u_1,...,\hat u_i,...,u_n$, let $E_i := \{\lambda_{i} = 0\}$ and $D_i = \{H_ i =0\}$ in $\mathbb C^n$. Note that for points in $E_{i,j} \cap \mathbb{C}_{I}$, we use $\lambda_i,H_{i},u_1,...,\hat u_i,...,\hat{u_{j}},...,u_n$ to be the local coordinates. Then by \textbf{Lemma \ref{slightly_reduction}} and \textbf{Theorem \ref{Thm_Formula_of_Zeta_Function_via_Embedded_Q_Resolution}}, $\mathbb L^{n}\cdot Z_{f}^{\mathrm{mot}}(s)$ is given as below. Note that $E_i \cap A_i = \varnothing$, we may drop some terms.
		\begin{align*}
			& \sum_{i=1}^n [A_i\setminus D_i] + \sum_{i=1}^n[A_i\cap D_i] \cdot \frac{(\mathbb L-1)\mathbb L^{-(s+1)}}{1-\mathbb L^{-(s+1)}} \\
			& \quad + \sum_{i=1}^n\sum_{I \subseteq Z_i} [(\mathbb C_I\setminus D_i)/ G_{i}] \cdot S_{G_{I,i}'}(\bm N_3, \bm \nu_3,s) \cdot \frac{(\mathbb L-1)\mathbb L^{-(ds+\vert \bm w\vert)}}{1-\mathbb L^{-(ds+\vert \bm w\vert)}}\\
			& \quad + \sum_{i=1}^n\sum_{I \subseteq Z_i} [(\mathbb C_I\cap E_{i,j}) / G_{i}] \cdot S_{G_{I,i,j}}(\bm N_4, \bm \nu_4,s) \cdot \frac{(\mathbb L-1)\mathbb L^{-(s+1)}}{1-\mathbb L^{-(s+1)}}\cdot \frac{(\mathbb L-1)\mathbb L^{-(ds+\vert \bm w\vert)}}{1-\mathbb L^{-(ds+\vert \bm w\vert)}}\\
			& = \sum_{i=1}^n [A_i\setminus D_i] + \sum_{i=1}^n[A_i\cap D_i] \cdot \frac{(\mathbb L-1)\mathbb L^{-(s+1)}}{1-\mathbb L^{-(s+1)}} \\
			& \quad + \sum_{i=1}^n\sum_{I \subseteq Z_i} ([\mathbb C_I/G_{i}]-[\hat E_{I,i} / G_{i}]) \cdot S_{G_{I,i}'}(\bm N_3, \bm \nu_3,s) \cdot \frac{(\mathbb L-1)\mathbb L^{-(ds+\vert \bm w\vert)}}{1-\mathbb L^{-(ds+\vert \bm w\vert)}}\\
			& \quad + \sum_{i=1}^n\sum_{I \subseteq Z_i} [(\mathbb C_I\cap E_{i,j}) / G_{i}] \cdot S_{G_{I,i,j}}(\bm N_4, \bm \nu_4,s) \cdot \frac{(\mathbb L-1)\mathbb L^{-(s+1)}}{1-\mathbb L^{-(s+1)}}\cdot \frac{(\mathbb L-1)\mathbb L^{-(ds+\vert \bm w\vert)}}{1-\mathbb L^{-(ds+\vert \bm w\vert)}}.
		\end{align*}
		Since $\Pi: \widehat{\mathbb C^n}(\bm w)\setminus \hat E \overset{\simeq}{\longrightarrow} \mathbb C^n\setminus\{\bm 0\}$, we see $\sum_{i=1}^n [(A_i\setminus D_i)/G_i] = \mathbb L^n - [Z(f)]$ and $\sum_{i=1}^n[(A_i\cap D_i)/G_i] = [Z(f)] - 1$. By \textbf{Lemma \ref{Quotient_Singularity_Grothendieck_Class}}, $[\mathbb C_I/G_{i}] = (\mathbb L-1)^{\vert I\vert}$ and we get $Z_{f}^{\mathrm{mot}}(s)$.
		
		For the local case, we only need to drop the strata that do not intersect $\hat E$. For the last assertion, one should simply notice that $\hat E_{I,i} = \{H_i = 0\}\cap \{\lambda_{i} = 0\} = \{f_d(\bm U_i)+\lambda_{i} h(\bm U_i^x) = 0\} \cap \{\lambda_{i} = 0\} = \{f_d(\bm U_i) = 0\} \cap \{\lambda_{i} = 0\}$. Hence the corresponding ``$\hat E_{I,i}$'' for $f_d$ coincides with that of $f$. 
	\end{proof}
	
	\subsection{Local Motivic Zeta Function in Dimension Two}
	
	In dimension $2$, we can write the local formula of \textbf{Theorem \ref{motivic zeta function}} very briefly.
	\begin{corollary}\label{dimension 2 result}
		Suppose $f \in \mathbb C[x_1,x_2]$ is semi-quasihomogeneous and has a unique singularity at $\bm 0$. Let $\bm w = (p,q)$, $(p,q) = 1$, be the weight of $f$ and $d = \deg_{\bm w} f$. Suppose $f(x_1^{p},x_1^{q}u) = x_1^d f_1(u)$ and $f(x_2^{p}v,x_2^{q}) =x_2^{d}f_2(v)$. Then
		\begin{align*}
			&\mathbb L^2\cdot Z_{f,\bm 0}^{\mathrm{mot}}(s) \\
			& = \frac{(\mathbb L-1)\mathbb L^{-(ds+\vert \bm w\vert)}}{1-\mathbb L^{-(ds+\vert \bm w\vert)}}\cdot 
			\begin{cases}
				\mathbb L-1-\frac{d}{pq}+S_1+S_2 + \frac{d}{pq}\frac{(\mathbb L-1)\mathbb L^{-(s+1)}}{1-\mathbb L^{-(s+1)}}, & f_1(0)\neq 0, f_2(0) \neq 0.\\
				\mathbb L-1-\frac{d-p}{pq}+S_2 + (\frac{d-p}{pq}+S_1')\frac{(\mathbb L-1)\mathbb L^{-(s+1)}}{1-\mathbb L^{-(s+1)}}, & f_1(0) = 0, f_2(0) \neq 0.\\
				\mathbb L-1-\frac{d-q}{pq}+S_1 + (\frac{d-q}{pq}+S_2')\frac{(\mathbb L-1)\mathbb L^{-(s+1)}}{1-\mathbb L^{-(s+1)}}, & f_1(0) \neq 0, f_2(0) = 0.\\
				\mathbb L-1-\frac{d-p-q}{pq}+(\frac{d-p-q}{pq}+S_1'+S_2')\frac{(\mathbb L-1)\mathbb L^{-(s+1)}}{1-\mathbb L^{-(s+1)}}, & f_1(0) = 0, f_2(0) = 0.
			\end{cases}	
		\end{align*}
		Here $S_i = S_{G_i}(\bm N_3,\bm \nu_3,s),\ S_i' = S_{G_i'}(\bm N_4,\bm \nu_4,s)$, and $G_1 = \frac{1}{p}(-1,q),\ G_1 = \frac{1}{q}(-1,p),\ G_1' = \frac{1}{p}(-1,d),\ G_2' = \frac{1}{q}(-1,d).$
		
		In particular, $Z_{f,\bm 0}^{\mathrm{mot}}(s)$  depends  only on the weight type and $\bm w$-degree of $f$.
	\end{corollary}
	\begin{proof}
		By \textbf{Theorem \ref{motivic zeta function}}, we may assume $f$ is quasihomogeneous and only compute $\hat E_{I,i}$.	
		
		We start with the cases that neither $x_1$ nor $x_2$ are factors of $f$. Then $f$ is convenient, i.e. $x_1^{a_1}$ and $x_2^{a_2}$ appears in $f$ for some $a_1,a_2 > 0$. In this case, it is equivalently that $x_1^{d/p},x_2^{d/q}$ appears in $f$. $f_1(u)$ is polynomial with non-zero constant term and of degree $p$. Similarly, $f_2(v)$ is of degree $d/p$ with $f_2(0) \neq 0$. Then $\hat{E}_{\varnothing,1} = \hat{E}_{\varnothing,2} = \varnothing$. Looking at \textbf{Theorem \ref{motivic zeta function}}, we only have to show $\hat E_{\{2\},1}/G_{I,1}$ consists of $\frac{d}{pq}$ points. Note that $E_{\{2\},1}/G_{I,i} \simeq \{f_1 = 0\}/\langle \zeta_p\rangle \subseteq \mathbb C^*/\langle \zeta_p\rangle$. $\langle \zeta_p \rangle$ acts on the roots of $f_1$ by multiplication, and hence faithfully. So it suffices to prove $f_1$ has no multiple roots. It is true since otherwise $f$ will have a multiple factor and is not an isolated singularity at $\bm 0$.
		
		For the second case, suppose $f = x_1 \tilde f$, with $\tilde f$ convenient. Then $\hat E_{\varnothing,1} = \{pt\}$ and $\hat{E}_{\varnothing,2} = \varnothing$. By the same argument as above, $\hat E_{\{2\},1}/G_{I,1}$ consists of $\frac{d-p}{pq}$, since only $\tilde f$ contribute to $\hat E_{\{2\},1}$. A similar result holds for the case $f = x_2\tilde f$, with $\tilde f$ convenient.
		
		For the last case, suppose $f = x_1x_2 \tilde f$, with $\tilde f$ convenient. Then $\hat E_{\varnothing,1} = \hat E_{\varnothing,2} = \{pt\}$ and we use the same argument again.
		
		As for the last assertion, one should notice that the classification of the four cases corresponds to (1) $pq\mid d$; (2) $p\mid d, q\nmid d$; (3) $q\mid d, p\nmid d$; (4) $p,q\nmid d$. So the classification is also given by the weight type and the $\bm w$-degree.
	\end{proof}
	
	From the formula, we are able to prove the survival of poles in dimension $2$.
	\begin{corollary}\label{dimension 2 non-cancellation}
		Suppose $f \in \mathbb C[x_1,x_2]$ is semi-quasihomogeneous and has a unique singularity at $\bm 0$. Let $\bm w = (p,q)$, $(p,q) = 1$, be the weight of $f$ and $d$ be the $\bm w$-degree of its principal part. Suppose $d > p+q$. Then $-1$ and $=\vert \bm w\vert/d$ are both poles of $Z_{f,\bm 0}^{\mathrm{mot}}(s)$.
	\end{corollary}
	\begin{proof}
		Since $d \geq pq > p+q = \vert \bm w\vert$, $-\vert \bm w\vert/d > -1$ is the negative of the log canonical threshold of $f$ at $\bm 0$. It is a pole of $Z_{f,\bm 0}(s)$ (see \cite{NX16}). As for $-1$, one simply notices that $\frac{d}{pq},\frac{d-p}{pq}+S_1',\frac{d-q}{pq}+S_2',\frac{d-p-q}{pq}+S_1'+S_2'$ are congruent to the sum of a combination of powers of $L$ with positive coefficient, modulo $\mathbb L^{s+1}$, respectively.
	\end{proof}
	\begin{remark}
		If $d \leq p+q$, then after a coordinate change, the principal part of $f$ is $x_1$ or $x_1x_2$. The result for these cases are trivial.
	\end{remark}
	\begin{remark}
		In a previous research about motivic principal values integral (MPVI, see \cite{Vanishing_of_Principal_Value_Integrals_on_Surfaces}) for hyperplane arrangements (\cite{MPVIHA}), it was shown that for generic hyperplane arrangement, the MPVI is always non-zero. Using the residue interpretation, we could also see $-2/d$ is a pole $Z_{f,\bm 0}(s)$ when $f \in \mathbb C[x_1,x_2]$ is a reduced generic hyperplane arrangement. Hence, the two results coincide in this case.
	\end{remark}

	\subsection{Examples and Discussion}
	
	\begin{example}
		Let $f = x^2+y^3 \subseteq \mathbb C[x,y]$. We have the following.
		
		\noindent(1) $H_1 = 1+y^3,H_2 = x^2+1$, $\bm w = (3,2)$. $G_{\{2\},1} = \{1\}$, $G_{\varnothing,1} = \mathrm{diag}\{\zeta_{3}^{-1},\zeta_3^2\}$, and $\mathrm{diag}\{(-1)^3,-1\}$. (2) $[Z(f)] = \mathbb L$. (3) $[(\mathbb C_{\varnothing,1}\setminus D_1)/G] = [(\mathbb C_{\varnothing,2}\setminus D_2)/G] = 1$ (only one point) and $[(\mathbb C_{\varnothing,1}\cap D_1)/G] = [(\mathbb C_{\varnothing,2}\cap D_2)/G] = 0$. (4) $[(\mathbb C_{\{2\}}\setminus D_1)/G] = \mathbb L-2$ and $[(\mathbb C_{\{2\}}\cap D_1)/G] = 1$. 
		
		Therefore, we have the zeta function of $f$.
		\begin{align*}
			& \mathbb L^2\cdot Z_{f}^{\mathrm{mot}}(s) \\
			& = \mathbb L(\mathbb L-1) + (\mathbb L-1) \frac{\mathbb L-1}{\mathbb L^{s+1}-1}+((\mathbb L-2)+(1+\mathbb L^{2s+2}+\mathbb L^{4s+4})+(1+\mathbb L^{3s+3}))\frac{\mathbb L-1}{\mathbb L^{6s+5}-1}\\
			& \quad + \frac{\mathbb L-1}{\mathbb L^{s+1}-1}\cdot \frac{\mathbb L-1}{\mathbb L^{6s+5}-1}\\
			& = \mathbb L(\mathbb L-1) \frac{\mathbb L^{6s+5}}{\mathbb L^{6s+5}-1} + (\mathbb L-1)^2 \frac{\mathbb L^{6s+5}}{(\mathbb L^{s+1}-1)(\mathbb L^{6s+5}-1)} + (\mathbb L^{2s+2}+\mathbb L^{4s+4}+\mathbb L^{3s+3}) \frac{\mathbb L-1}{\mathbb L^{6s+5}-1}\\
			& =  (\mathbb L-1) \frac{\mathbb L\cdot \mathbb L^{6s+5}\cdot \mathbb L^{s+1} -\mathbb L^{6s+5}}{(\mathbb L^{s+1}-1)(\mathbb L^{6s+5}-1)} + (\mathbb L^{2s+2}+\mathbb L^{4s+4}+\mathbb L^{3s+3}) \frac{\mathbb L-1}{\mathbb L^{6s+5}-1}.
		\end{align*}
		It coincides with the one in \cite{An_Introduction_to_p_adic_and_Motivic_Integration_Vius_Sos}.
	\end{example}
	
	\begin{example}
		Let $f = x^5+yz^2+xy^3 \subseteq \mathbb C[x,y,z]$. 
		
		The weight of $f$ is $\bm w = (6,8,11)$ and $d = 30$. We have $H_1 = 1+yz^2+y^3$, $H_2 = x^5+z^2+x$, $H_3 = x^5+y+xy^3$. We shall compute each bracket above.
		
		\begin{table}[h!]
			\centering
			\begin{tabular}{|c|c|c|c|}
				\hline
				$i$ & $I$ & $[(\mathbb C_I\cap D_i)/G_{I,i}]$ & $G_{I,i}$\\
				\hline
				$1$ & $\{2,3\}$ & $\mathbb L-3$ & $\{1\}$\\
				\hline
				$1$ & $\{2\}$ & $1$ ($1+y^3$, three points contract to $1$) & $\frac{1}{2}(1,0,1)$\\
				\hline
				$1$ & $\{3\}$ & $0$ ($1$, no root) & $\{1\}$\\
				\hline
				$1$ & $\varnothing$ & $0$ & $\frac{1}{6}(-1,8,11)$\\
				\hline
				$2$ & $\{3\}$ & $1$ ($z$, only one point) & $\{1\}$\\
				\hline 
				$2$ & $\varnothing$ & $1$ & $\frac{1}{8}(6,-1,11)$\\
				\hline
				$3$ & $\varnothing$ & $1$ & $\frac{1}{11}(6,8,-1)$\\
				\hline
			\end{tabular}
		\end{table}
		For $I = \{2,3\}$, we identify $(\mathbb C^*)^2/\frac{1}{6}(8,11)$ with $(\mathbb C^*)^2/\frac{1}{6}(1,-1)$ via
		\begin{displaymath}
			[(u,v)] \mapsto [(uv^2,v)],\ [(s,t)] \mapsto [(t^{-2}s,t)].
		\end{displaymath} 
		Under this isomorphism, $f$ is sent to $1+s+t^{-6}s^3$ in fake coordinate $s,t$ on $(\mathbb C^*)^2/\frac{1}{6}(1,-1)$. $s,t^6$ is the coordinate of $(\mathbb C^*)^2 \simeq (\mathbb C^*)^2/\frac{1}{6}(1,-1)$. Hence $[(C_I\cap D_i)/G_i] = [(\mathbb C^*)^2 \cap \{a^2+a^2b+b^3 = 0\}]$. Next we solve the equation $a^2+a^2b+b^3 = 0$, when $a,b \neq 0$. Under isomorphism $a\mapsto ab$ and $b\mapsto b$, $[(\mathbb C^*)^2 \cap \{a^2+a^2b+b^3 = 0\}] = [(\mathbb C^*)^2 \cap \{a^2+a^2b+b = 0\}]$. Hence, if $a\neq \pm \sqrt{-1}$, $b$ is killed. Besides, there is no root when $a = \pm 1$. Therefore, the class is equal to $\mathbb L-3$. Then $Z_{f}^{\mathrm{mot}}$, $Z_{f,\bm 0}^{\mathrm{mot}}$ can be given accordingly.

	\end{example}
	
	\begin{example} \label{homogeneous motivic zeta function}
		Let $f \in \mathbb C[\bm x]\setminus \mathbb C$ be a homogeneous polynomial of degree $d$, with a unique singularity at the origin. Then $(\mathbb C^n,f)$ is (log) resolved by blowing up at the origin. By \textbf{Theorem \ref{Ch2_Sec3_Motivic_and_Log_Resolution_Thm}}, we have
		\begin{displaymath}
			Z_{f,\bm 0}^{\mathrm{top}}(s) = \frac{1}{s+1}\bigg(\chi(\mathbb P^{n-1})-\chi(Z_+(f)) +   \frac{\chi(Z_+(f)}{ds+n}\bigg) = Z_{f,\bm 0}^{\mathrm{top}}(s) = \frac{1}{s+1}\bigg(n-\chi(Z_+(f)) +  \frac{\chi(Z_+(f)}{ds+n}\bigg).
		\end{displaymath}
		Where $Z_+(f) = \{[\bm a] \in \mathbb P^{n-1} \mid f(\bm a) = 0\}$. Since $Z_+(f)$ is smooth, apply Proposition 10.4.1 in \cite{Intersection_homology_perverse}, we see
		\begin{displaymath}
			\chi(Z_+(f)) = n -\frac{1}{d}\bigg(1+(-1)^{n-1}(d-1)^{n} \bigg),
		\end{displaymath}
		which is determined by $d$.
	\end{example}
	
	Based on  \textbf{Example \ref{homogeneous motivic zeta function}} and \textbf{Corollary \ref{dimension 2 result}}, we propose the  \noindent\textbf{Conjecture \ref{topological zeta function topological}} (it is included in the introduction section).
	
	
	\section{Stringy E-function via Embedded $\mathbb{Q}$-Resolution}

	\subsection{Formula of Stringy E-function}\label{sec7}
	
	The techniques of embedded $\mathbb Q$-resolutions can be applied to stringy E-function. 
	\begin{theorem}\label{stringy E-function for Q-Gor}
		Assume $X$ is a $\mathbb{Q}$-Gorenstein variety with at worst log terminal singularity of pure dimension $n$. Let $\pi : Y \longrightarrow X$ be an embedded $\mathbb{Q}$-resolution such that $Y=\bigsqcup_{k \ge 0} Y_{k}$, and for any point $q \in Y_{k}$, we have $(Y,q) \simeq (\mathbb{C}^{n}/G_{k}, 0)$, and the relative canonical divisor $K_{\pi}$ is locally given by $x_{1}^{\nu_{1,k}-1}...x_{n}^{\nu_{n,k}-1}$, where $x_{1},...,x_{n}$ are the coordinates of $\mathbb{C}^{n}$. Then the stringy E-function of $X$ is given by:
		\begin{flalign*}
			E_{st}(X)=\sum_{k \ge 0}E(Y_{k})E(S_{G_{k}}(\bm{\nu}_{k}))\prod_{i=1}^{n}\frac{uv-1}{(uv)^{\nu_{i,k}}-1}.
		\end{flalign*}
		\begin{proof}
			Take a log resolution of $X$, $\phi: Y' \longrightarrow X$. By the formula of \cite{Introduction_to_Motivic_Integration}, we can rewrite the definition of stringy E-function into the integration as follows:
			\begin{flalign*}
				E_{st}(X)=E\big( \mathbb{L}^{n}\int_{J_{\infty}(Y')}\mathbb{L}^{-\mathrm{ord}_{t}K_{\phi}}\,\mathrm{d}\mu_{J_{\infty}(Y')}\big).
			\end{flalign*}
			Since $Y$ is smooth, the $\mathbb{Q}$-Gorenstein measure and the usual measure on $J_{\infty}(Y')$ coincide. By the change of variables formula in \cite{Motivic_Zeta_Function_on_Q_Gorenstein_Varieties_Leon_Martin_Veys_Viu_Sos}, we get:
			\begin{flalign*}
				E_{st}(X)&=E(\mathbb{L}^{n}\int_{J_{\infty}(Y')}\mathbb{L}^{-\mathrm{ord}_{t}K_{\phi}}\,\mathrm{d}\mu_{J_{\infty}(Y')}^{\mathbb{Q}\mathrm{-Gor}}) = E(\mathbb{L}^{n}\int_{J_{\infty}(X)}1\,\mathrm{d}\mu_{J_{\infty}(X)}^{\mathbb{Q}\mathrm{-Gor}}).
			\end{flalign*}
			Since we have an embedded $\mathbb{Q}$-resolution $\pi: Y \longrightarrow X$, we can use the change of variables formula in \cite{Motivic_Zeta_Function_on_Q_Gorenstein_Varieties_Leon_Martin_Veys_Viu_Sos} again, and calculating the resulting integral by    in \cite{Motivic_Zeta_Function_on_Q_Gorenstein_Varieties_Leon_Martin_Veys_Viu_Sos}:
			\begin{flalign*}
				E_{st}(X)&=E(\mathbb{L}^{n}\int_{J_{\infty}(Y)}\mathbb{L}^{-\mathrm{ord}_{t}K_{\pi}}\,\mathrm{d}\mu_{J_{\infty}(Y)}^{\mathbb{Q}\mathrm{-Gor}}) =\sum_{k \ge 0}E(Y_{k})E(S_{G_{k}}(\bm{\nu}_{k}))\prod_{i=1}^{n}\frac{uv-1}{(uv)^{\nu_{i,k}}-1}.
			\end{flalign*}
		\end{proof}
	\end{theorem}
	
	\subsection{The Stringy E-Function for Semi-quasihomogeneous Polynomial}\label{sec8}
	Now we fix a semi-quasihomogeneous polynomial $f \in \mathbb{C}[x_{1},...,x_{n}]$ with a unique singularity at $0$ and we start to calculate the stringy E-function of $H=\{f=0\}$ using this formula. We use the notations in \textbf{Subsection \ref{subse_zeta_function}}. By \textbf{Theorem \ref{Thm_Formula_of_Zeta_Function_via_Embedded_Q_Resolution}}, $\pi$ is a $\mathbb Q$-resolution.
	
	Let $\hat{E}:=E \cap \hat{H}$ be the exceptional divisor of $\pi$. We consider the points $q$ in the chart $U_{i}-\bigcup_{j=1}^{i-1}U_{j}$. When $q \notin \hat{E}$, we have $(\hat{H},q) \cong (\mathbb{C}^{n-1},0)$ and the corresponding $\nu_{i}=1$ for all $i$, so the contribution of these point to the sum is:
	\begin{flalign*}
		E(U_{i}\cap \hat{H}-\bigcup_{j=1}^{i-1}(U_{j} \cap \hat{H})).
	\end{flalign*}
	We denote the part of $(U_{i}-\bigcup_{l=1}^{i-1}U_{l}) \cap \hat{E}$ where $\frac{\partial f_{d}}{\partial u_{j}} \neq 0,\frac{\partial f_{d}}{\partial u_{l}}=0$($j,l \neq i, 1 \le l \le j-1$) by $E_{i,j}$.
	
	When $q \in \hat{E}$, there exists $j$ such that $q \in E_{i,j}$, and we have 
	$(\hat{H},q) \cong (\mathbb{C}^{n-1}/G_{i,I,j},0)$ for some $G_{i,I,j}$ if $q \in \mathbb{C}_{I}$. The contribution to the sum is:
	\begin{flalign*}
		\sum_{j \neq i}\sum_{I \subset \{1,...,\hat{i},...,\hat{j},...,n\} }E((E_{i,j} \cap \mathbb{C}_{I})/G_{i})E(S_{G_{i,I,j}}(\bm{|w|},1,...,1))\frac{uv-1}{(uv)^{\bm{|w|}}-1}.
	\end{flalign*}
	
	If we take the sum of all contributions from the part where $q \notin \hat{E}$, we have:
	\begin{flalign*}
		&\sum_{i=1}^{n}E(U_{i}\cap \hat{H}-\bigcup_{j=1}^{i-1}(U_{j} \cap \hat{H}))  \\
		=&E(\hat{H}-\hat{E}) \\
		=&E(H)-1.
	\end{flalign*}
	Finally, we have the following theorem:
	\begin{theorem}\label{stringy E-function for semi-quasihomogeneous}
		The stringy invariant of $H$ is given by:
		\begin{flalign*}
			E(H)-1+\sum_{i=1}^{n}\sum_{j \neq i}\sum_{I \subset \{1,...,\hat{i},...,\hat{j},...,n\} }E(E_{i,j} \cap \mathbb{C}_{I})E(S_{G_{i,I,j}}(\bm{|w|},1,...,1))\frac{uv-1}{(uv)^{\bm{|w|}}-1}.
		\end{flalign*}
	\end{theorem}

	\subsection{Newton Polyhedron, Non-degenerate Polynomials and Normal Fan}
	
	For stringy E-function for non-degenerate polynomials, authors in \cite{Stringy_E-Function_of_Hypersurface} have given a formula in section $3$. But we will give a brand-new formula for it. In this subsection, we recall some basic facts about non-degenerate polynomials, Newton polyhedron and its normal fan. 
	
	We fix a polynomial $f \in \mathbb{C}[x_{1},...,x_{n}]$, and for $\bm{a}=(a_{1},...,a_{n}) \in \mathbb{Z}_{\ge 0}^{n}$, let $\bm{x}^{\bm{a}}=x_{1}^{a_{1}}...x_{n}^{a_{n}}$ so we can write $f$ as $f=\sum_{\bm{a}}c_{\bm{a}}\bm{x}^{\bm{a}}$. We define $Supp(f):=\{\bm{a}|c_{\bm{a}} \neq 0\}$, and the Newton polyhedron of $f$, denoted by $\Gamma(f)$, is defined to be the convex hull of $\bigcup_{\bm{a} \in Supp(f)}\bm{a}+\mathbb{R}^{n}_{\ge 0}$. A face $\tau$ of $\Gamma(f)$ is a nonempty set of the form $\tau=H\cap \Gamma(f)$, where $H$ is a hyperplane of $\mathbb{R}^{n}$ such that $\Gamma(f)$ is on the one side the halfspace determined by $H$. For every face $\tau$ of $\Gamma(f)$, we define $f_{\tau}:=\sum_{\bm{a} \in \tau}c_{\bm{a}}\bm{x}^{\bm{a}}$, and we say $f$ is non-degenerate if $f_{\tau}$ is non-singular in $(\mathbb{C}-\{0\})^{n}$ for every face $\tau$ of $\Gamma(f)$. Now we associate a $\mathbb{R}$-linear function $\phi_{f}$ to $\Gamma(f)$ as follows:
	\begin{flalign*}
		\phi_{f}: \mathbb{R}_{+}^{n}& \longrightarrow \mathbb{R}_{\ge 0}  \\
		\bm{v}&  \mapsto inf_{\bm{a} \in \Gamma(f)}\bm{v} \cdot \bm{a}.
	\end{flalign*}
	Using $\phi_{f}$, we define the normal fan of $\Gamma(f)$, denoted by $\Sigma(f)$, as follows:
	\begin{flalign*}
		\Sigma(f):=\{\sigma_{\bm{a}}|\;\bm{a} \in \mathbb{R}_{+}^{n}\},
	\end{flalign*}
	where
	\begin{flalign*}
		\sigma_{\bm{a}}:=\{\bm{v} \in \mathbb{R}^{n}_{+}|\; \phi_{f}(\bm{v})=\bm{a} \cdot \bm{v}\}.
	\end{flalign*}
	Note that $\Sigma(f)$ is a subdivision of $\mathbb{R}^{n}_{+}$, so this subdivision induces a morphism between toric varieties associated to $\Sigma(f)$ and $\mathbb{R}^{n}_{+}$: $X_{\Sigma(f)} \longrightarrow X_{\mathbb{R}^{n}_{+}}=\mathbb{C}^{n}$.

		\begin{center}
			\begin{tikzpicture}
				\draw [->](0,0)--(3.5,0);
				\node[above] at (3.5,0) {x};
				\draw [->](0,0)--(0,3);
				\node[right] at (0,3) {y};
				\draw (0,2.5)--(1,1);
				\draw (1,1)--(3,0);
				\node[below] at (3,0) {(6,0)};
				\node[right] at (0,2.5) {(0,5)};
				\node[above right] at (1,1) {(2,2)};
				\draw[dashed][->] (0,0)--(0.5,1);
				\draw[dashed][->] (0,0)--(0.75,0.5);
			\end{tikzpicture}
			\hspace{-3cm} Figure 1.
	\end{center}

	\begin{example}
We consider the non-degenerate polynomial $f=x^6+x^{2}y^{2}+y^{5}$. The Newton polyhedron and its normal fan are drawn in the following picture, where the two dashed lines are the vertical lines of the two one-dimensional faces of $\Gamma(f)$. The normal fan $\Sigma(f)$ is actually the subdivision of $\mathbb{R}^{2}_{+}$ by these two dashed lines (see Figure 1).
		\end{example}

	\subsection{Stringy E-Function for Non-degenerate Polynomial}\label{sec10}
	We fix a non-degenerate polynomial $f \in \mathbb{C}[x_{1},...,x_{n}]$ which defines a unique singularity at $0$ with at most log terminal singularity, let $\Gamma(f)$ be its Newton polyhedron and $\Sigma(f)$ be its normal fan. We choose a simplicial subdivision of $\Sigma(f)$ without changing the set of 1-dim cones of it, denoted by $\Sigma(f)'$, and by Section 3.1.6 in \cite{Motivic_Zeta_Function_on_non-Degenerate_Polynomial}, the toric resolution $\Pi: X_{\Sigma(f)'} \longrightarrow X_{\mathbb{R}^{n}_{+}}=\mathbb{C}^{n}$ is an embedded $\mathbb{Q}$-resolution. We explain this resolution in details. Suppose $\sigma \subset \Sigma(f)'$ is a cone with dimension $n$, and since it is simplicial, we assume that it is generated by $\rho_{i}=(a_{i1},...,a_{in}) \in \mathbb{Z}_{\ge 0}^{n}$ for $1 \le i \le n$ with $gcd(a_{i1},...,a_{in})=1$. The coordinate ring of the affine variety $U_{\sigma}$ is given by $\mathbb{C}[\check{\sigma}\cap \mathbb{Z}^{n}]=\mathbb{C}[x_{1}^{b_{1}}...x_{n}^{b_{n}}|b_{1}a_{i1}+...+b_{n}a_{in} \ge 0, \forall i]$. Let $c_{i}=b_{1}a_{i1}+...+b_{n}a_{in}, \forall i$, and this is equivalent to a coordinate change $u_{1}^{c_{1}}...u_{n}^{c_{n}}=x_{1}^{b_{1}}...x_{n}^{b_{n}}$. So the condition about $b_{i}$ is transformed into the condition of $c_{i}$ such that $c_{i} \in \mathbb{Z}_{\ge 0}$ and they satisfy some congruence conditions about $A$ to assure that $b_{i} \in \mathbb{Z}$, where $A=det(a_{ij})$. This tells us that $U_{\sigma} \simeq \mathbb{C}^{n}/G_{\sigma}$ for some finite abelian group $G_{\sigma}$ with coordinates $u_{1},...,u_{n}$. In these new coordinates, $x_{i}=u_{1}^{a_{1i}}...u_{n}^{a_{ni}}$, and we can rewrite $f$ in the new coordinates: for $\bm{a} \in Supp(f)$, $c_{\bm{a}}\bm{x}^{\bm{a}}=c_{\bm{a}}u_{1}^{\bm{a}\cdot \rho_{1}}...u_{n}^{\bm{a}\cdot \rho_{n}}$. Consider the following set $\sigma^{*} \subset \Gamma(f)$:
	\begin{flalign*}
		\sigma^{*}:=\{\bm{a} \in \Gamma(f)| \;\forall \bm{v} \in \sigma, \phi_{f}(\bm{v})=\bm{a} \cdot \bm{v}\}.
	\end{flalign*}
	Note that $\sigma^{*}$ is a face of $\Gamma(f)$. Since $\rho_{i} \in \sigma$, $m_{i}:=inf_{\bm{a} \in Supp(f)}\bm{a} \cdot \rho_{i}$ is taken by those $\bm{a} \in \sigma^{*}$. Now we assume $f=u_{1}^{m_{1}}...u_{n}^{m_{n}}g_{\sigma}$, where $g_{\sigma}:=\sum_{\bm{a} \in \sigma^{*}}c_{\bm{a}}+\sum_{\bm{a} \notin \sigma^{*}}c_{\bm{a}}u_{1}^{\bm{a}\cdot \rho_{1}-m_{1}}...u_{n}^{\bm{a}\cdot \rho_{n}-m_{n}}$. Section 3.1.2 in \cite{Motivic_Zeta_Function_on_non-Degenerate_Polynomial} shows that the strict transform of $f$ in $U_{\sigma}$ is given by $g_{\sigma}$, the exceptional part of the resolution is given by $u_{1}^{m_{1}}...u_{n}^{m_{n}}$, and $g_{\sigma}$ is smooth in $\mathbb{C}^{n}$.
	
	Now we are able to give the formula of stringy E-function of $H=\{f=0\} \subset \mathbb{C}^{n}$. We use the embedded $\mathbb{Q}$-resolution $\Pi: X_{\Sigma(f)'} \longrightarrow \mathbb{C}^{n}$. We denote the strict transform of $H$ by $\hat{H}$ and $\pi:=\Pi|_{\hat{H}}:\hat{H} \longrightarrow H$. Suppose $\Sigma(f)'=\bigcup_{i=1}^{s}\sigma_{s}$ and $\sigma_{i}$ is a cone of dimension $n$, so $X_{\Sigma(f)'}=\bigcup_{i=1}^{s}U_{\sigma_{i}}=\bigsqcup_{i=1}^{s}(U_{\sigma_{i}}-U_{\sigma_{1}}\cup...\cup U_{\sigma_{i-1}})$. Since $g_{\sigma_{i}}$ is smooth , $\pi$ is also an embedded $\mathbb{Q}$-resolution of $H$. Similar to the semi-quasihomogeneous case, we introduce the notations of the stratification. Let $V_{i}:=U_{\sigma_{i}}-U_{\sigma_{1}}\cup...\cup U_{\sigma_{i-1}}$ be a closed subset of $U_{\sigma_{i}}$, and $E_{i,j}$ be the subset of $V_{i} \cap \hat{H}$ such that $\frac{\partial g_{\sigma_{i}}}{\partial u_{j}} \neq 0$ and $\frac{\partial g_{\sigma_{i}}}{\partial u_{1}}=...=\frac{\partial g_{\sigma_{i}}}{\partial u_{j-1}}=0$. For $I \subset \{1,...,\hat{j},...,n\}$ and point $p \in \mathbb{C}_{I}\cap E_{i,j}$, we have $(\hat{H},p)\simeq (\mathbb{C}^{n-1}/G_{i,I,j},0)$ for some $G_{i,I,j}$. Next we compute the restriction of $K_{\pi}$ on $U_{\sigma_{i}}$. Assume the $1$-dim cones of $\sigma_{i}$ is given by $\rho_{i,1},...,\rho_{i,n}$, and the coordinates of $U_{\sigma_{i}}$ are given by $u_{1},...,u_{n}$ as the last paragraph. The irreducible components of $K_{\pi}$ on $U_{\sigma_{i}}$ are given by $u_{1}=0,...,u_{n}=0$, and the coefficients are $\rho_{i,l} \cdot (1,1,...,1)-\phi_{f}(\rho_{i,l})-1$(the calculation of these coefficients are the same as the proof of Proposition 2.3 in \cite{Stringy_E-Function_of_Hypersurface}). 
	
	In conclusion, we have the following formula for stringy E-function of $H$:
	\begin{theorem}\label{stringy E-function for nondegenerate}
		The stringy E-function of $H$ is given by:
		\begin{flalign*}
			\sum_{i=1}^{s}\sum_{j=1}^{n}\sum_{I \subset \{1,...\hat{j},...,n\}}E((E_{i,j}\cap \mathbb{C}_{I})/G_{\sigma_{i}})S_{G_{i,I,j}}(\bm{\nu}_{i,I,j})\prod_{l \in I}\frac{uv-1}{(uv)^{\rho_{i,l}\cdot (1,...,1)-\phi_{f}(\rho_{i,l})}-1},
		\end{flalign*}
		where the $l$-th coordinate of $\bm{\nu}_{i,I,j}$ is $\rho_{l}\cdot (1,...,1)-\phi_{f}(\rho_{l})$ if $l \in I$, and $1$ otherwise.
	\end{theorem}
	
	\begin{remark}
		By section 4.6 in \cite{Stringy_E-Function_of_Non-Degenerate_Polynomial}, the log discrepancy of $H$ are integers, so the log terminal condition is equivalent to canonical. Thus one can apply Proposition $2.3$ in \cite{Stringy_E-Function_of_Hypersurface} and the calculation of motivic zeta function of $H$ in \cite{Motivic_Zeta_Function_on_non-Degenerate_Polynomial} to get the stringy E-function of $H$.
	\end{remark}

	\end{document}